\documentclass{article}

\usepackage[T1]{fontenc}
\usepackage[utf8]{inputenc}
\usepackage{lmodern}
\usepackage{microtype}
\usepackage[margin=1.in]{geometry}
\usepackage{amsmath}
\usepackage{amssymb}
\usepackage{amsthm}
\usepackage{mathtools}
\usepackage{bbm}
\usepackage{enumerate}
\usepackage{enumitem}
\usepackage[dvipsnames]{xcolor}

\definecolor{blue}{HTML}{2980b9}
\definecolor{gray}{HTML}{bdc3c7}
\definecolor{purple}{HTML}{c0392b}
\definecolor{green}{HTML}{2ecc71}
\definecolor{yellow}{HTML}{f1c40f}

\usepackage{multirow}
\usepackage{empheq}
\usepackage{listings}
\usepackage{color}
\usepackage[%
    font={small,sf},
    labelfont=bf,
    format=hang,    
    format=plain,
    margin=0pt,
    width=0.8\textwidth,
]{caption}
\usepackage{subcaption}
\usepackage{comment}

\usepackage[square,comma,numbers]{natbib}
\usepackage{hyperref}
\usepackage{graphicx} % Required for including pictures

\usepackage{tikz}
\usetikzlibrary{babel}

\numberwithin{equation}{section}
\usepackage{pgfplots}
\pgfplotsset{width=10cm,compat=1.9}

\usepackage{environ}
\makeatletter
\newsavebox{\measure@tikzpicture}
\NewEnviron{scaletikzpicturetowidth}[1]{%
  \def\tikz@width{#1}%
  \begin{lrbox}{\measure@tikzpicture}%
  \BODY
  \end{lrbox}%
  \pgfmathparse{#1/\wd\measure@tikzpicture}%
  \BODY

}
\makeatother

\usepackage[dvipsnames]{xcolor}
\usepackage{soul}
\usepackage{ upgreek }

\newcommand*\indica[1]{\mathbbm{1}_{ #1 }}
\def\cA{{\mathcal A}}

\def\cF{{\mathcal F}}
\def\cG{{\mathcal G}}

\def\cI{{\mathcal I}}

\def\cL{{\mathcal L}}
\def\cM{{\mathcal M}}
\def\cN{{\mathcal N}}

\def\cP{{\mathcal P}}

\def\cX{{\mathcal X}}

\def\C{\mathbb{C}}

\def\E{\mathbb{E}}
\def\F{\mathbb{F}}

\def\P{\mathbb{P}}
\def\R{\mathbb{R}}

\def\T{\mathbb{T}}

\def\RR{\R}

\def\d{\mathrm{d}}
\def\balpha{\boldsymbol{\alpha}}

\def\mfS{\mathfrak{S}}

\providecommand{\abs}[1]{\lvert#1\rvert}
\providecommand{\norm}[1]{\lVert#1\rVert}

\def\unig{m \otimes g}

\usepackage{MnSymbol}
\theoremstyle{plain}
\newtheorem{theorem}{Theorem}[section]
\newtheorem{lemma}[theorem]{Lemma}
\newtheorem{corollary}[theorem]{Corollary}
\newtheorem{proposition}[theorem]{Proposition}

\newtheorem{definition}[theorem]{Definition}
\newtheorem{example}[theorem]{Example}
\newtheorem{remark}[theorem]{Remark}

\usepackage[symbol]{footmisc}

\usepackage{fancyhdr}
\date{\vspace{-1em}\normalsize{\today}}

\title{Kuramoto Mean Field Game with Intrinsic Frequencies}
\author
{Rene Carmona\footnote{Research partially supported by AFOSR under Grant No.FA9550-23-1-0324.
Department of Operations Research and Financial
Engineering, Princeton University, Princeton, NJ, 08540, USA, email: 
{\tt rcarmona@princeton.edu}}
\and Quentin Cormier\footnote{CMAP, Ecole Polytechnique,
Palaiseau, Paris,  France, email: 
{\tt quentin.cormier@inria.fr}}
\and H. Mete Soner \footnote{Research partially supported by the National Science Foundation grant DMS 2406762. 
Department of Operations Research and Financial
Engineering, Princeton University, Princeton, NJ, 08540, USA, email: 
{\tt soner@princeton.edu}.}}
\date{\today}
\begin{document}
\maketitle
\abstract{
	This paper studies a mean field game formulation of the classical
	Kuramoto model for synchronization. Our model
	captures the diversity within the population by considering
	random intrinsic frequencies, which allows us to study
	the impact of this heterogeneity on synchronization patterns and stability. Our findings contribute insights into the interplay between intrinsic frequency diversity and synchronization dynamics, offering a more realistic understanding of complex systems. The proposed framework has broad applications ranging from coupled oscillators in physics to social dynamics, and serves as a valuable tool for studying networks with distributed intrinsic frequencies.
}
\vspace{10pt}
\smallskip\newline
\noindent\textbf{Key words:} Mean field games, Kuramoto model, Synchronization,
viscosity solutions.
%\vspace{10pt}
\smallskip\newline
\noindent\textbf{Mathematics Subject Classification:}  35Q89, 35D40, 39N80, 91A16, 92B25
\vspace{10pt}
%\tableofcontents

%%%%%%%%%%%%%%%%%%%%%%%%%%
\section{Introduction}
\label{se:intro}

This paper uses mean field games
\cite{Car,CD,huang_large_2006,huang_large-population_2007,huang_nash_2007,LL1,LL2,LL3}
to study synchronization as initially modelled in \cite{YMMS,YMMS2}.
We refer the reader to our earlier work \cite{CCS}
for historical background on mean field games,
the role of the classical Kuramoto
model in the understanding of spontaneous synchronization phenomena,
as well as the introduction and the first mathematical analysis of the model
recast as a Mean Field Game (MFG) for a family of oscillators in mean field interaction
with homogenous frequencies.

Models that encourage synchronization typically go through a phase transition
from incoherent states to self-organization as
the incentive for concentration gets larger~\cite{CCS, CesaroniCirant}. This fascinating phenomenon
is the key feature of the classical Kuramoto model and its generalization
to the game theoretic approach is the central focus of this study.
In the mean field Kuramoto game,
random movement of the oscillators or players
on the one-dimensional torus $\T$ are driven
by their intrinsic frequencies drawn from a given distribution $g$ on the real line
and the distribution of the others.  The scaling parameter $\kappa$
determines the relative strength of these opposing
effects, with large $\kappa$ values achieving coordination.
{One of the main contributions of the present paper is to show the existence of phase transitions by identifying specific values of this parameter delineating regimes with different equilibrium behavior.}
They are subject to a Brownian noise of strength $\sigma$, and
control their drift so as to minimize a cost functional
defined on an infinite horizon with a discount factor $\beta>0$.

We aim to develop a comprehensive
analysis of Kuramoto synchronization games
by studying the existence of non-trivial stationary solutions
and the stability of the incoherent state.
We construct several critical thresholds
that govern these results which are analogous to the ones
identified  in the classical
dynamical systems setting.
In addition to its intrinsic value,
Kuramoto games also serve as a
sophisticated example for the study of
the long time behavior of mean field games.
While there are several deep results for Lasry-Lions monotone games
or small interactions ~\cite{cecchin2024exponentialturnpikephenomenonmean},
without these assumptions one expects exciting structures with
multiple invariant distributions, periodic solutions or even
more complicated attractors \cite{Ci}.
These issues also arise in a similar
one-dimensional model for congestion with a local
interaction studied in \cite{gomes2016,gomes2018}. There, the authors consider deterministic optimal control
problems with ergodic cost. Additionally, a class of problems with
trigonometric running cost is studied  in \cite{n2017} where some explicit solutions
are constructed.

Recently, \cite{MR4053591} and \cite{MR3936909}  provide a novel approach
to potential MFGs without monotonicity or small interactions
by employing the celebrated weak KAM theory.
In addition, \cite{BC} and \cite{CH} provide deep insights into the large-time behavior of potential MFGs, including the convergence of the fictitious play algorithm.
Here, we rely on an alternative approach to the stability of stationary solutions by exploiting the implicit function theorem and a fine analysis of the action of the main operator in the Laplace/Fourier domain.

In the context of the Kuramoto model, incoherence manifests itself by
oscillators having the uniform distribution $m$ on the torus.
Indeed, it is shown in Lemma~\ref{lem:uniform-measure} below
that the product measure $m\otimes g$  is always
a Nash equilibrium, {or equivalently a solution 
of the MFG. Here, and throughout the paper, we use the notation
$m$ for the uniform measure on the torus, namely $m(\d x)= \d x/(2\pi)$, and $g$ is the notation we use for the distribution of the natural frequencies of the oscillators.} At a critical
threshold, this incoherent solution loses its stability
and other equilibria emerge.
%\qc{This is a minor detail but here I would say simply instead: ``other solutions emerge''; because when the incoherent solution losses its stability, then periodic solutions may appear instead of stationary solutions, that is we can have $\kappa_1(g) > \kappa_2(g)$.  }
As the non-trivial spatial structure
of these solutions indicates synchronization among the particles,
the phase transition
is related to {either the appearance
	of non-uniform stationary equilibria as defined in Definition \ref{de:MFG_solution},
	or to the local stability of the equilibrium
	$\unig$ as defined  in Definition \ref{def:stable}. The intuitive meaning of the product probability measure $m\otimes g$ is that the state frequencies $x$ are uniformly distributed while independently, the intrinsic frequencies $\omega$ have distribution $g$.
	Accordingly, these asymptotic regimes are identified by the critical values of the parameter $\kappa$ defined by:}
\begin{align}
	\label{eq:kappa1}
	\kappa_1(g) & :=
	\inf\{ \kappa\ :\  \ \text{There are non-uniform
		                     stationary Nash equilibria at}\
	\kappa',\ \ \forall \kappa' \ge \kappa\}, \\
	\label{eq:kappa2}
	\kappa_2(g) & := \sup\{\kappa\ :\
	\unig\ \text{is stable at}\ \kappa'\
	\text{in the sense of Definition \ref{def:stable},}\ \forall \kappa' \le \kappa\}.
\end{align}

Our first main results, Theorems \ref{th:ub} and \ref{th:lb2} below, imply that for
frequency distributions with a positive Fourier transform,
there exists a \emph{unique} {value of the parameter $\kappa$, say $\kappa_c(g)$,} satisfying,
\begin{equation}
	\label{fo:tilde_kappa_c}
	\kappa_1(g)\le
	\kappa_c(g) \le\kappa_2(g),
	\quad \text{where}\quad
	\kappa_c (g):=  \left( \int_{\R }
	\frac{  \gamma \sigma^2 + 2 \omega^2 }{(\gamma^2 + \omega^2)(\sigma^4 + 4 \omega^2)}
	\ g(\d \omega)\ \right)^{-1},
\end{equation}
and $\gamma :=\beta+ \sigma^2/2$.
These inequalities suggest a phase
transition at the critical parameter $\kappa_c(g)$.

A simpler model with two states studied in \cite{HS}
allows for a complete analysis;
the uniform distribution is the only stationary
solution for small values and it is globally stable,
while non-uniform stable equilibria emerge
in the supercritical case.  Moreover, the analogous
critical parameters are all equal to each other.
However, the stability structure is
far richer for general distributions
and analysis on the Laplace domain provides a
deeper understanding of the system.
We further identify in Section~\ref{sec:Penrose} the important \emph{Penrose critical threshold} $\kappa_P(g)$
given by,
\begin{equation}
	\label{eq:ckappa-penrose}
	\kappa_P(g) :=
	\inf \{ \kappa > 0\ : \ \exists \theta \in \R
	\ \text{such that}\ P(i \theta) = 2 /\kappa \},
\end{equation}
where
\begin{equation}
	\label{eq:P}
	P(z) := \int_{\R} \frac{1}{(\gamma+i \omega - z)(\frac{\sigma^2}{2} + z - i \omega )}
	\ g(\d \omega),\qquad \text{on} \ \mfS
	:= \{ z \in \C: \quad
	\Re(z) \in  (-\frac{\sigma^2}{2}, \gamma) \}.
\end{equation}
Here, $\sigma, \gamma > 0$ are parameters of the model
introduced earlier, and the above function $P$ is related to the Laplace
transform of the linearization of the solution
map around $\unig$ derived in Section~\ref{sec:Frechet}.

Our next main result, Theorem~\ref{th:two_Diracs},
proves that $\kappa_2(g) \ge \kappa_P(g)$
when $g$ is the sum of two Dirac measures.
Note that for any symmetric distribution $g$, it holds that $P(0) = 2 / \kappa_c(g)$
and so
\[ \kappa_P(g) \le\kappa_c(g). \]
Numerically, we observe that
the  equality holds when $g$ has a positive
Fourier transform, an equality which we can prove for $g=\delta_0$.
In general, these thresholds do not agree
as shown in Example \ref{ex:counter2}
when $g$ is the sum of two point masses
revealing the complexity of the solutions
of the Kuramoto games.

The paper is organized as follows.
After  introducing the model
and the main results in the next section,
we prove the results in the super-critical regime in
Section~\ref{sec:stationary},
and Section~\ref{s.stable}
investigates the stability of the uniform equilibrium
for the sub-critical case.
Section~\ref{sec:Penrose} outlines
an approach based on the Laplace transform,
and the case of two Dirac measures
is studied in Section~\ref{s.two_Diracs}.
Sections~\ref{sec:estimates} and~\ref{sec:Frechet}
provide the technical analysis for the
derivation of the Fr\'echet derivative which is central to our
stability results.
Related coupled partial differential equations and
the potential structure are given in  Appendix~\ref{app:pde}.
In the subsequent short appendices, we state technical results used in the paper.

%%%%%%%%%%%%%%%%%%%%%%%%%%
\section{Model and Results}
\label{se:The_problem}
The one-dimensional torus is denoted by $\T := \R / (2\pi \mathbb{Z})$.
Throughout the paper, we denote by $\cP(A)$ the set of probability measures on $A$. The probability measure $g \in \cP(\R)$ with a finite first moment
is the distribution of the intrinsic frequencies of the oscillators.
For $x \in \T$ and $\mu \in \cP(\T \times \R)$, we set
\begin{equation}
	\label{fo:c_of_mu}
	c(x, \mu) := \int_{\T \times \R} 2\sin^2(\frac{x-y}{2}) \mu(\d y, \d \omega).
\end{equation}
Note that this function only involves the first marginal of $\mu$ and
\begin{equation}
	\label{fo:developed_c}
	c(x, \mu) = 1 - \cos(x) \int_{\T} \cos(y) \mu(\d y, \R) - \sin(x) \int_{\T} \sin(y) \mu(\d y, \R).
\end{equation}

We fix a filtered probability space $(\Omega,\F,\P)$ supporting an $\F$-adapted Brownian motion $(B_t)_{t \geq 0}$, and we assume that the filtration $\F = (\cF_t)_{t \geq 0}$ satisfies the usual conditions, i.e., $\cF_0$ is complete and $\cF_t$ is right-continuous. We assume that the initial $\sigma$-field $\cF_0$ is such that for any probability measure $\nu_0 \in \cP(\T)$, one can construct an $\cF_0$-measurable, $\T$-valued random variable $X_0$ with distribution $\nu_0$.
We denote by $\cA_t$ the set of all square integrable $\F$-progressively measurable controls $\alpha: [t, \infty)\times \Omega \rightarrow \R$ and write $\cA := \cA_0$.

Let $(\mu_t)_{t \geq 0}$ be a flow of probability measures such that the second marginal of each $\mu_t \in \cP(\T \times \R)$
%\rc{Why $\R_+$?} \ms{no reason, mistake}
is  equal to $g$:
\[
	\quad \mu_t(\T, \d \omega) = g(\d \omega), \qquad \forall t \geq 0.
\]
Then, in view of the disintegration theorem, for $g$-almost every
$\omega \in \R$ there exists a probability measure $\mu^\omega_t \in \cP(\T)$ such that:
\begin{equation}
	\label{eq:decomposition-mu_t}
	\mu_t(\d x, \d \omega) = \mu^\omega_t(\d x)
	g(\d \omega), \qquad \forall t \geq 0.
\end{equation}
For each $\alpha \in \cA$ and $\omega \in \R$, we consider the cost
\[ J^{\omega, \mu}(\alpha) := \E \int_0^\infty e^{-\beta t} \left[ \frac{1}{2} \alpha^2_t + \kappa c(X^{\omega, \alpha}_t, \mu_t) \right] \d t, \]
where $\kappa \geq 0$ and $\beta > 0$ are parameters of the model, and $X^{\omega, \alpha}_t$ is given by
\[
	X^{\omega, \alpha}_t = X^{\omega}_0 + \int_0^t \alpha_s \d s + \omega t + \sigma B_t,
\]
with initial state $X_0^\omega$ having distribution $\mu_0^\omega$, i.e., $\cL(X^\omega_0) = \mu^\omega_0$.
With these notations in hand, we define the solutions of the Mean Field Game (MFG) problem as follows:

\begin{definition}
	\label{de:MFG_solution}
	{\rm{We say that $(\mu_t)_{t \geq 0} \subset \cP(\T \times \R)$ is}}
    	a solution of the MFG problem {(or equivalently a Nash equilibrium)}
    with intrinsic frequency distribution $g$,
	{\rm{if the following hold:
			\begin{itemize}
				\item for all $t\ge 0$ the second marginal of $\mu_t$  is $g$, i.e., $\mu_t(\T, \cdot) =g$;
				\item  there exists a Borel set $E \subset \R$ with $g(E) = 1$
				      such that  for all $\omega \in E$ the disintegration formula \eqref{eq:decomposition-mu_t} holds
				      with measures $\mu^\omega_t \in \cP(\T)$,
				      and there exists $\alpha^\omega_* \in \cA$ such that
				      $$
					      \inf_{\alpha \in \cA} J^{\omega, \mu}(\alpha) = J^{\omega, \mu}(\alpha^\omega_*),
					      \qquad \text{and} \qquad
					      \cL(X^{ \omega, \alpha^\omega_*}_t) = \mu^\omega_t,\qquad t\ge 0.
				      $$
			\end{itemize}}}
\end{definition}
The stationary MFG equilibria are defined as follows:
\begin{definition}
	{\rm{A probability measure $\mu \in \cP(\T \times \R)$ is a}} stationary MFG
	equilibrium (or equivalently a Nash equilibrium)
    with intrinsic frequency distribution $g$, {\rm{if the constant flow $(\mu_t)_{t \geq 0}$ defined by
	$\mu_t = \mu$ for all $t\ge 0$  is a solution of the MFG problem in the sense of Definition \ref{de:MFG_solution} above.}}
\end{definition}
As is common in the literature, an alternate definition through
coupled Hamilton-Jacobi-Bellman and Kolmogorov-Fokker-Planck equations
is also available.  We discuss these approaches and the potential
structure of the game in Appendix \ref{app:pde}. Although our analysis
does not utilize these important properties, they provide valuable
insights.

\begin{remark}
	\label{rem:coupling0}
	The coupling between intrinsic frequencies in the above definition is subtle.
	It is through the interaction
	term $c(X^{\omega,\alpha}, \mu_t)$ which considers the distribution of all
	frequencies $\mu_t \in \cP(\T \times \R)$  and not only the ones with a fixed $\omega$. We should also note that the conditional probability
	flow $\mu^\omega_t$
	is not necessarily a Nash equilibrium for the mean-field game with $\omega$ fixed.
	Analogous coupling exists in the coupled equations and
	the potential problem introduced in Appendix~\ref{app:pde}, and
	further discussed in Remark~\ref{rem:coupling}.
\end{remark}
{
\begin{remark}
\label{rem:integrated_formulation}
An equivalent formulation of Definition~\ref{de:MFG_solution} can be established by treating the intrinsic frequency as a random variable $\Omega$. Let $(X_0, \Omega)$ be drawn from the joint initial distribution $\mu_0$, and define the enlarged filtration $\cG_t := \sigma(X_0, \Omega, B_s, s \leq t)$. Denoting by $\cA_{\cG}$ the space of $\cG_t$-progressively measurable controls, the state dynamics are given by $X_t = X_0 + \int_0^t \alpha_s \d s + \Omega t + \sigma B_t$. A flow $(\mu_t)_{t \geq 0}$ is then a solution if the second marginal is $g$, and there exists a global optimal control $\alpha_* \in \cA_{\cG}$ minimizing the cost
\[
    J(\alpha) := \E \int_0^\infty e^{-\beta t} \left[ \frac{1}{2} \alpha^2_t + \kappa c(X_t, \mu_t) \right] \d t,
\]
such that the joint law satisfies $\cL(X_t, \Omega) = \mu_t$ for all $t \geq 0$.
\end{remark}
}
In the simpler model with no intrinsic frequencies studied in~\cite{CCS},
the uniform distribution is always a stationary solution.  The analogue
of this distribution is the following
$g$-uniform measure $m\otimes g \in \cP(\T \times \R)$ given by,
\begin{equation}
	\label{eq:uniform}
	(m\otimes g)(\d x, \d \omega):= m(\d x) g(\d \omega),
	\qquad \text{where} \qquad
	m(\d x):=
	\frac{1}{2\pi}\ \d x.
\end{equation}
We then have the following immediate result.

\begin{lemma}
	\label{lem:uniform-measure}
	For all $\kappa \geq 0$ and $\beta, \sigma > 0$,
	the $g$-uniform distribution $\unig$
	is a stationary solution of the MFG problem 
    {(or equivalently a Nash equilibrium)}
    with intrinsic frequency distribution $g$.
\end{lemma}
\begin{proof}
	It is clear that $\unig$ satisfies
	the disintegration formula \eqref{eq:decomposition-mu_t}
	with $\mu^\omega_t=m$
	and  $\kappa c(x, \unig) = \kappa$ is independent of $x$.
	Therefore, $J^{\omega, \unig}(\alpha) \geq J^{\omega, \unig}(0)$
	for all $\alpha \in \cA$.
	Additionally,
	$$
		\cL(X^{\omega, 0}_t) = \cL(X^\omega_0 + \omega t +  \sigma B_t) =m,
	$$
	where we used that  $\cL(X^\omega_0) = m$.
	This shows that $\unig$ is a solution  of the MFG problem.
\end{proof}

\subsection{Rotation and Galilean invariance}
\label{ss.invaraince}

Let $(\uptheta, \tau) \in \T \times \R$ and consider the following maps
$$
	R_{\uptheta, \tau}(x, \omega) := (x + \uptheta, \omega + \tau), \quad r_\uptheta(x) := x + \uptheta, \quad t_\tau(\omega) = \omega+\tau,
$$
where  on the torus,
the addition is always mod $2\pi$.
Given $\mu_t \in \cP(\T \times \R)$ and $g \in \cP(\R)$, let $\tilde{\mu}_t$ and $\tilde{g}$
be the push forward probability measures of $\mu_t$ by $R_{\uptheta + \tau t, \tau}$ and of $g$ by $t_\tau$ respectively.
In other words:
\[ \tilde{\mu}_t :=  {R_{\uptheta + \tau t, \tau}}_{\#} \mu_t, \quad \tilde{g} := {t_\tau}_\# g. \]
Then we have:
\begin{lemma}
	\label{lem:symmetries}
	Let $(\mu_t)_{t\ge 0}$ be a solution of the MFG problem 
    {(or equivalently a Nash equilibrium)}
    with intrinsic frequency distribution $g$.
	Then $(\tilde{\mu}_t)$ is a solution of the MFG problem with intrinsic frequency distribution $\tilde{g}$.
\end{lemma}
\begin{proof}
	As $\mu_t = \mu^\omega_t(\d x) g(\d \omega)$,
	$$
		\tilde{\mu}_t = \tilde{\mu}^{\tilde{\omega}}_t(\d \tilde{x}) \tilde{g}(\d \tilde{\omega}) \quad \text{ with }
		\quad \tilde{\mu}^{\tilde{\omega}}_t := {r_{\uptheta + \tau t}}_{\#}  \mu^{\tilde{\omega} - \tau }_t.
	$$
	We set $\tilde{E} := \{ \omega+\tau, \omega \in E \}$,
	$\tilde{\omega} := \omega+ \tau \in \tilde{E}$, so that
	$$
		c(X^{\omega, \alpha}_t + \uptheta + \tau t, \tilde{\mu}_t) = c(X^{\omega, \alpha}_t, \mu_t),
	$$
	and consequently, $J^{\omega+\tau, \tilde{\mu}}(\alpha) = J^{\omega, \mu}(\alpha)$,
	for every $\alpha \in \cA$.
	These imply that
	$$
		\inf_{\alpha \in \cA} J^{\tilde{w}, \tilde{\mu}}(\alpha) = J^{\omega, \mu}(\alpha^\omega_*) = J^{\tilde{\omega}, \tilde{\mu}}(\alpha^{\tilde{\omega}-\tau}_*).
	$$
	Moreover, by definition,
	$$
		\cL(X^{\omega, \alpha^\omega_*}_t + \uptheta + \tau t)
		=  {r_{\uptheta + \tau t}}_{\#} \mu^\omega_t= \tilde{\mu}^{\tilde{\omega}}_t,
	$$
	proving that $(\tilde{\mu}_t)$ is a solution of the MFG problem
	with intrinsic frequency distribution $\tilde{g}$.
\end{proof}

\begin{remark}
	\label{re:rotation}
	By choosing $\tau = 0$, we find that any rotation by $\uptheta$ of a solution of the MFG problem
	is again a solution of the MFG problem with the same intrinsic frequency distribution.
	So even for a fixed intrinsic frequency distribution, stationary solutions of the MFG problem are only
	determined up to rotations on the torus.
	In addition, by choosing $\tau = -\int_{\R} \omega g(\d \omega)$, we see that $\tilde{g}$ has mean zero.
	So from now on, we assume without any loss of generality that this translation is always implemented,
	and \emph{we implicitly assume that $g$ is centered}, namely:
	\[ \int \omega g(\d \omega) = 0. \]
\end{remark}

\subsection{Main Results}
\label{ss.results}

Recall the critical thresholds $\kappa_1(g)$ of \eqref{eq:kappa1},
$\kappa_2(g)$ of \eqref{eq:kappa2},  $\kappa_c(g)$ of \eqref{fo:tilde_kappa_c},
$\kappa_P(g)$ of \eqref{eq:ckappa-penrose}, and the constant
\begin{equation}
	\label{eq:def-gamma}
	\gamma :=\beta+ \sigma^2/2\ ;
\end{equation}
all defined in the Introduction.
By definition,
the uniform measure $m \otimes g$ is the
unique solution when $\kappa <\kappa_1(g)$
and there may be multiple solutions for $\kappa  \ge \kappa_1(g)$.

In the classical Kuramoto model, larger values of $\kappa$
support coherence, and this general rule generalizes
to our game theoretic approach as well.
Indeed, when there are no intrinsic frequencies
or equivalently for $g(\d \omega) = \delta_0(\d \omega)$,
\cite{CCS} identifies $\kappa_c(\delta_0) = \gamma \sigma^2$
as the unique critical interaction strength of the model.
In particular, the  results of \cite{CCS} can be summarized as
$$
	\kappa_1(\delta_0)\le \kappa_c(\delta_0) \le\kappa_2(\delta_0).
$$
This paper is concerned with bounds of this type.
Our first result, proved in subsection~\ref{ss.proof_ub} and stated below, shows that
the critical threshold $\kappa_c(g)$ is an upper bound for $\kappa_1(g)$.
\begin{theorem}[Upper Bound - general]
	\label{th:ub}
	For any symmetric distribution $g$, {the uniform measure $m \otimes g$ is the
unique solution for any $\kappa < \kappa_c(g)$. Consequently,
$\kappa_1(g) \le \kappa_c(g)$.}
\end{theorem}
\vspace{5pt}

Our second set of results, established  in Section~\ref{s.stable},
provides lower bounds for $\kappa_2(g)$;
first statement below follows from Lemma~\ref{lem:first},
and the second from Lemma~\ref{lem:second}.
\begin{theorem}[Lower bound - general]
	\label{th:lb1}
	For any symmetric distribution $g$, $\kappa_2(g) \ge \kappa_c(\delta_0)=\gamma\sigma^2$.
\end{theorem}
\noindent
To state the next result, we recall that the
Fourier transform of $g$ is given by,
\begin{equation}
	\label{eq:def-fourier-transform-g}
	\mathfrak{F}g(t):= \hat{g}(i t),
	\qquad \text{where}\qquad
	\hat{g}(z) := \int_{\R} e^{-  z\omega} g(\d \omega),
	\qquad z \in \C.
\end{equation}
\begin{theorem}
	\label{th:lb2}
	Suppose that $\mathfrak{F}g(t) \ge 0$ for every $t \in \R$.  Then,
	$\kappa_2(g) \ge \kappa_c(g)$.
\end{theorem}

These bounds combined identify
$\kappa_c(g)$ as  \emph{the unique critical threshold}
of the Kuramoto game for  distributions with a positive
Fourier transform:
$$
	\kappa_1(g) \le \kappa_c (g)  \le \kappa_2(g),
	\qquad \text{for any} \ g \ \text{satisfying}\
	\mathfrak{F}g \ge0.
$$
In particular, the above inequalities represent a direct generalization
of the results of \cite{CCS}.
Additionally, numerically
we have observed that they hold as equalities
for distributions with positive Fourier transforms.

Our second main theorem considers the case when $g$ is concentrated at finite
frequencies. Its proof is presented in Subsection~\ref{ss.temporal} after Proposition~\ref{prop:inversibility-I-L}. The analysis in
Section~\ref{s.two_Diracs} is related to the Penrose constant and employs
Laplace transform techniques.
{
\begin{theorem}
	\label{th:two_Diracs}
	Suppose that $g(\d \omega)= \sum_{i = 1}^N \alpha_i (\delta_{\omega_i} + \delta_{-\omega_i})$,
	for some $N \geq 1$, $\omega_i \geq 0$ and $\alpha_i > 0$,
    with $\sum_{i = 1}^N \alpha_i=1/2$.  Then, $\kappa_2(g) \ge \kappa_P(g)$.
\end{theorem}
}

\begin{example}
	\label{ex:positive}
	{\rm{Any Gaussian  with mean zero satisfies the hypothesis of
			Theorem~\ref{th:lb2}.  With $\sigma= 2$, $\beta= 1$ and the
			variance of Gaussian equal to $1$, we have
			$$
				\kappa_c (g)   = 10 \sqrt{\frac{2}{\pi }} \left(e^{9/2} \text{erfc}\left(\frac{3}{\sqrt{2}}\right) +e^2\text{erfc}\left(\sqrt{2}\right)\right)^{-1}
				\approx 13.77 > \kappa_c(\delta_0) = 12.
			$$
			In Example \ref{ex:counter2} below, we compute the Penrose
			constant $\kappa_P(g)$ and $\kappa_c (g)$ when $g$ is as in Theorem~\ref{th:two_Diracs}.
		}}
\end{example}

%%%%%%%%%%%%%%%%%%%%%%%%%%%%%%
\section{Stationary solutions}
\label{sec:stationary}
In the subsequent subsections, we first characterize
the stationary solutions through a two-dimensional fixed point
problem and then
prove Theorem~\ref{th:ub}.

%%%%%%%%%%%%%%%%%%%%
\subsection{Description of the stationary solutions}
\label{ss.describe_stationary}
The central observation  of our analysis
is the characterization of the stationary solutions of the MFG problem
by two real numbers $\alpha_1$ and  $\alpha_2$.
For each coupling constant $\kappa>0$ and
$\balpha=(\alpha_1,\alpha_2)\in\RR\times\RR$, we define the function $c(\cdot,\balpha)$ by:
\begin{equation}
	\label{fo:c_alpha}
	\kappa c(x, \balpha) := \kappa - \alpha_1 \cos(x) - \alpha_2 \sin(x).
\end{equation}
The above notation is consistent with~\eqref{fo:developed_c}
as long as $\balpha=(\alpha_1,\alpha_2)$ is derived from the probability measure $\mu$ as follows:
\begin{equation}
	\label{fo:alphas}
	\alpha_1 := \kappa \int \cos(y) \mu(\d y, \R)   \quad \text { and } \quad \alpha_2 := \kappa \int \sin(y) \mu(\d y, \R).
\end{equation}
For each $\omega\in\R$, we denote by $v^\omega$
the unique classical solution of the following stationary HJB equation:
\begin{equation}
	\label{eq:HJB-omega-eq}
	\omega \partial_x v^\omega + \frac{\sigma^2}{2} \partial_{xx} v^\omega + \kappa c(x, \balpha) - \frac{1}{2} (\partial_x v^\omega)^2
	= \beta v^\omega.
\end{equation}
If $\alpha_1$ and $\alpha_2$ are given by \eqref{fo:alphas}, this function gives the optimal value of a generic player with intrinsic frequency $\omega$ responding to the environment (field) of the other players as given by the measure $\mu$.
%In equilibrium or equivalently when the players are following the optimal strategy, 
The optimal  dynamics of the state of this generic player are given by the following equation:
\begin{equation}
	\label{eq:optimal-state-player-w}
	X^\omega_t = X^\omega_0 + \omega t - \int_0^t  \partial_x v^\omega(X^\omega_s) \d s + \sigma B_t.
\end{equation}
Note that the right hand side of the above equation is understood modulo $2\pi$, that is, as an element of
$\mathbb{T}=\mathbb{R}/(2\pi\mathbb{Z})$.
The following is proved in \cite{CG}.
\begin{lemma}
	\label{lem:stationary-solution-torus}
	If $v^\omega$ is the unique solution of the HJB equation \eqref{eq:HJB-omega-eq}, the stochastic differential equation~\eqref{eq:optimal-state-player-w} has a unique invariant probability measure on $\T$ given by
	\begin{equation}
		\label{fo:nu_omega}
		\nu^\omega_\infty(x) \d x = \frac{1}{Z \xi^\omega(x)} \left[ \int_0^{2 \pi} \xi^\omega(y) \d y  + (e^{-\frac{4 \pi}{\sigma^2} \omega } - 1) \int_0^x \xi^\omega(y) \d y \right] \d x,
	\end{equation}
	where $Z$ is the normalizing factor and
	\begin{equation}
		\label{eq:def-of-xi-omega}
		\xi^\omega(x) := e^{-\frac{2}{\sigma^2}(\omega x - v^\omega(x))}.
	\end{equation}
\end{lemma}

% Note that $v^\omega$ and $\nu^\omega_\infty$ only depend on $\mu_\infty$ thought the values of $\alpha_1$ and $\alpha_2$, which appear in the cost $c(x, \mu_\infty)$. 
Here and in the following we use the same notation $\nu^\omega_\infty$ for the measure and its density.
For each value $\kappa$ of the coupling parameter, we define the function $F_\kappa$ from $\RR\times\RR$ into itself by:
\begin{equation}
	\label{fo:F_kappa}
	F_{\kappa}(\alpha_1, \alpha_2) :=  \kappa \left( \int_{\T \times \R} \cos(y) \nu^\omega_\infty(y)  g(\d \omega) \d y, \int_{\T \times \R} \sin(y) \nu^\omega_\infty(y)  g(\d \omega) \d y  \right)
\end{equation}
where $\nu^\omega_\infty$ is the invariant probability measure \eqref{fo:nu_omega} with $v^\omega$ solving \eqref{eq:HJB-omega-eq} for $\balpha=(\alpha_1,\alpha_2)$. These allow us to state the characterization result.

\begin{proposition}[Characterization of stationary equilibria]
	\label{pr:characterization}
	There is a one-to-one correspondence between the fixed points of $F_\kappa$ and the stationary solutions of the MFG problem.
\end{proposition}
\begin{proof}
	Let $(\alpha_1, \alpha_2)$ be a fixed point of $F_\kappa$. We claim that
	\begin{equation}
		\label{fo:mu_infinity}
		\mu_\infty(\d x, \d \omega) = \nu^\omega_\infty(x) \d x \;g (\d \omega)
	\end{equation}
	is a stationary solution of the MFG problem. Indeed, for all $\omega \in \R$, $\alpha^\omega_*(t) = -\partial_x v^\omega(X^\omega_t)$ is an optimal control for a generic player with intrinsic frequency $\omega$, since $v^\omega$ is the smooth solution of \eqref{eq:HJB-omega-eq} associated to $\balpha$ satisfying \eqref{fo:alphas} with $\mu(\d y, \d \omega)=\nu^\omega_\infty(\d y)g(\d \omega)$, because $\balpha=(\alpha_1,\alpha_2)$ is a fixed point of $F$. In addition, if  $\cL(X^\omega_0) = \nu^\omega_\infty$, the solution of \eqref{eq:optimal-state-player-w} has law $\cL(X^\omega_t) = \nu^\omega_\infty$, proving that $\mu_\infty$ defined in \eqref{fo:mu_infinity} is a stationary MFG equilibrium.
	\vskip 2pt
	Conversely, if $\mu$ is a stationary MFG equilibrium with intrinsic frequency distribution $g$, defining $c(\cdot,\mu)$ by~\eqref{fo:c_of_mu} and using the corresponding $\alpha_1$ and $\alpha_2$, so that $c(x,\mu)=c(x,\balpha)$, for each $\omega\in\T$, the classical solution $v^\omega$ of the stationary HJB equation~\eqref{eq:HJB-omega-eq} is the equilibrium value function giving the optimal feedback control function $\alpha^*(x)=-\partial_xv^\omega(x)$ leading to the dynamics given by the SDE~\eqref{eq:optimal-state-player-w} whose unique invariant distribution is $\nu^\omega_\infty(x)dx$. This implies that
	\[
		\mu_\infty(\d x, \d \omega) = \nu^\omega_\infty(\d x) \;g (\d \omega)
	\]
	and hence, that $(\alpha_1,\alpha_2)$ is a fixed point for the function $F_\kappa$ defined in~\eqref{fo:F_kappa}.
\end{proof}

%%%%%%%%%%%%%%%%%%
\subsection{Linearization}

In the previous subsection, we reduced the solutions to stationary
MFGs to two-dimensional fixed point problems.
Additionally, in Lemma~\ref{lem:uniform-measure} we have
shown that the  $g$-uniform distribution $\unig$ is a stationary
solution corresponding to the fixed point $(0, 0)$ of $F$.
These observations allow us to construct other solutions by
linearization around the origin.

We start this analysis
by computing the derivative of $F_\kappa$ at the origin.
To achieve this, we linearize the  equation~\eqref{eq:HJB-omega-eq} around
the particular solution $v=\kappa/\beta$ corresponding to the $g$-uniform equilibrium $\unig$.
The corresponding linear equation is derived
by first considering the equation satisfied by $v^\omega(x)-\kappa/\beta$
and then dropping the quadratic term. The result is:
$$
	\omega \partial_x u^\omega + \frac{\sigma^2}{2} \partial_{xx} u^\omega = \beta u^\omega + \ell(x),
	\qquad\text{where}\qquad
	\ell(x):=\alpha_1 \cos(x) + \alpha_2 \sin(x).
$$
By the Feynman-Kac formula,
$$
	u^\omega(x) = -\E \int_0^\infty e^{-\beta s} \ell(X^\omega_s) \d s, \qquad \text{ with }
	\qquad X^\omega_s = x + \omega s + \sigma B_s.
$$
Recall $\gamma$ of \eqref{eq:def-gamma}.
Using the explicit expression of $\ell(x)$, we find that
$$
	u^\omega(x) = - \frac{\alpha_1 \gamma + \alpha_2 \omega}{\gamma^2 + \omega^2} \cos(x) + \frac{\alpha_1 \omega - \alpha_2 \gamma}{\gamma^2 + \omega^2} \sin(x) =: A \cos(x) + B \sin(x).
$$
We now use Lemma~\ref{lem:stationary-solution-torus} to arrive at:

\begin{lemma}
	For $\kappa>0$ fixed, as $\alpha_1, \alpha_2$ tend to zero,
	$$
		\nu^\omega_\infty(x) = \frac{1}{2 \pi} + \frac{1}{2 \pi} \frac{\lambda}{(1+\lambda^2)\omega} \left[ (B \lambda - A)\cos(x) - (B + \lambda A) \sin(x) \right] + o(\abs{\alpha_1} + \abs{\alpha_2}),
	$$
	where $\lambda := 2 \omega/\sigma^2$ and $o(\cdot)$ is
	the classical notation representing a function satisfying
	$\lim_{r \downarrow 0} o(r)/|r| =0$.
\end{lemma}
\begin{proof}
	We have
	\begin{align*}
		\xi^\omega(x) = e^{-\lambda x}(1 + \frac{\lambda}{\omega} (A \cos(x) + B \sin(x))) + o(\abs{A} + \abs{B}).
	\end{align*}
	Therefore, we find that
	\begin{align*}
		\nu^\omega_\infty(x) \propto \frac{1 - e^{-2 \pi \lambda}}{\lambda} + \frac{(1-e^{-2 \pi \lambda})(B \lambda - A)}{(1+\lambda^2)\omega} \cos(x) - \frac{(1-e^{-2 \pi \lambda})(B + \lambda A)}{(1+\lambda^2)\omega} \sin(x)
		+ o(\abs{A} + \abs{B}).
	\end{align*}
	The normalizing factor is $Z = \frac{1- e^{-2 \pi \lambda}}{\lambda} 2 \pi$.
	We conclude the proof by combining the above statements.
	The justification that the correction terms are in $o(\abs{\alpha_1} + \abs{\alpha_2})$ is done exactly as in \cite[Lemma 7.1]{CCS}.
\end{proof}
We leverage the above asymptotic expression to compute the derivative of
$F_\kappa$ at the origin.
\vspace{10pt}

\begin{proposition}
	\label{prop:one-to-one-stat-solution-2d-function}
	The function $F_\kappa$ is
	differentiable at $(0,0)$ and
	$$
		\d F_\kappa(0,0) \cdot \begin{pmatrix} \alpha_1 \\ \alpha_2 \end{pmatrix}
		=  \int_{\R } \frac{\kappa}{(\gamma^2 + \omega^2)(\sigma^4 + 4 \omega^2)}
		\begin{pmatrix} \gamma \sigma^2 + 2 \omega^2 & (\sigma^2 - 2 \gamma)\omega
		                \\ (2 \gamma - \sigma^2)\omega & \gamma \sigma^2 + 2 \omega^2\end{pmatrix} \begin{pmatrix}\alpha_1 \\ \alpha_2
		\end{pmatrix}
		\  g(\d \omega)
	$$
	where the $\cdot$ in the left hand side stands for matrix multiplication.
\end{proposition}
\vspace{10pt}

%%%%%%%%%%%%%%%%%%%%%%%%%%%%%
\subsection{Symmetric distributions}
\label{ss.symmetric}
We now assume that $g$ is symmetric around zero, that is $g(-\d \omega) = g(\d \omega)$.
We first reduce the problem to a scalar fixed-point problem.
\begin{lemma}
	Assume that $\alpha_2 = 0$ such that the cost is $c(x, \mu_\infty) = \kappa -\alpha_1 \cos(x)$, $\alpha_1 \in \R$. The value function solving \eqref{eq:HJB-omega-eq}, and the invariant probability measure given by Lemma~\ref{lem:stationary-solution-torus} satisfy for all $x \in \T$
	\[ v^\omega(-x) = v^{-\omega}(x), \quad \text{ and } \quad \nu^\omega_\infty(-x) = \nu^{-\omega}_\infty(x). \]
\end{lemma}
\begin{proof}
	Using the fact that $c(x, \mu_\infty)$ is even, we find that $v^\omega(-x)$ and $v^{-\omega}(x)$ satisfy the same HJB equation. By uniqueness, it follows that $v^\omega$ satisfies $v^\omega(-x) = v^{-\omega}(x)$, from which it follows that
	$\xi^\omega(-x) = \xi^{-\omega}(x)$. Recall that the expression of $\xi^\omega$ is given in~\eqref{eq:def-of-xi-omega}. Using the periodicity of $v^\omega$ we get $\xi^{-\omega}(z-2 \pi) = \xi^{-\omega}(z) e^{-\frac{4 \pi \omega}{\sigma^2}}$, and therefore:
	\begin{align*}
		\int_0^{2 \pi} \xi^\omega(y) \d y & = \int_0^{2 \pi} \xi^{-\omega}(-y) \d y =  \int_0^{2 \pi} \xi^{-\omega}(z - 2 \pi) \d z = e^{-\frac{4 \pi \omega}{\sigma^2}} \int_0^{2 \pi} \xi^{-\omega}(z) \d z .
	\end{align*}
	Finally, using the explicit expression of $\nu^\omega_\infty(x)$,
	we conclude that $\nu^\omega_\infty(-x) = \nu^{-\omega}_\infty(x)$ as claimed.
\end{proof}
For a symmetric $g$,
\[
	\int_{\T \times \R} \sin(y) \nu^\omega_{\infty}(y) g(\d \omega) = 0.
\]
Therefore, the function $F_\kappa$ maps
the line $\{(\alpha_1, 0);\, \alpha_1 \in \R \}$ to itself.
We set
\[ G_\kappa(\alpha) := F_\kappa(\alpha, 0)_1
	= \kappa \int_{\T \times \R} \cos(y) \nu^\omega_{\infty}(y) g(\d y). \]

\begin{proposition}[Characterization of stationary equilibria for symmetric $g$]
	Assume that $g$ is symmetric around zero.
	Then, there is a one-to-one correspondence between
	the fixed points of $G_\kappa: \R_+  \rightarrow \R_+$
	and the stationary solutions of the MFG problem.
	In addition, $G_\kappa$ is differentiable at $\alpha=0$ with
	$$
		G_\kappa'(0) = \frac{\kappa}{\kappa_c(g)},
		\qquad \text{where} \ \kappa_c(g)\ \text{is as in
			\eqref{fo:tilde_kappa_c}}.
	$$
\end{proposition}
\begin{proof}
	The form of $G_\kappa'(0)$ is a particular case of the formula obtained in
	Proposition~\ref{prop:one-to-one-stat-solution-2d-function}.
	If $\alpha$ is a fixed point of $G_\kappa$ then $(\alpha, 0)$ is a fixed point of $F_\kappa$ and so by Proposition~\ref{pr:characterization}, there exists an associated stationary solution of the MFG problem.
	Consider now $\mu$ to be a stationary solution of the MFG problem. Then using Lemma~\ref{lem:symmetries} with $\tau = 0$ and $\uptheta \in \R$ well-chosen,  we can rotate $\mu$ such that	\[ \int \cos(x) \mu(\d x, \d \omega) \geq 0, \quad  \int \sin(x) \mu(\d x, \d \omega) = 0. \]
	By definition, there exists $E \subset \R$ with $g(E) = 1$ and a family of probability measures $\mu^\omega_\infty(\d x), \omega \in E$ such that
	$\mu = \mu^\omega(\d x) g(\d \omega)$.
	For each $\omega \in E$, a player with intrinsic frequency $\omega$ has a unique optimal control given by $\partial_x v^\omega$ and one has
	\[ \cL(X^\omega_t) = \mu^\omega = \nu^\omega_\infty, \]
	as $\nu^\omega_\infty$, given by Lemma~\ref{lem:stationary-solution-torus}, is the unique invariant distribution of \eqref{eq:optimal-state-player-w}.

	Therefore, $\alpha = \kappa \int \cos \mu(\d x, \d \omega) \geq 0$ is a fixed point of $G_\kappa$.
\end{proof}

\subsection{Proof of Theorem~\ref{th:ub}}
\label{ss.proof_ub}
For $\kappa > \kappa_c$,
by the explicit formula we conclude that $G_\kappa'(0) > 1$.
As $G_\kappa$ is continuous and $G_\kappa \leq \kappa$,
this implies that $G_\kappa$ has at least a fixed point different from zero.

\qed

\begin{example}
	{\rm{If $g(\d \omega) =  \indica{[-a, a] }(\omega) \frac{\d \omega}{2a}$
			is the uniform measure on $[-a, a]$ for some $a > 0$, then}}
	\[ \kappa_c(g) = \frac{\arctan \left(\frac{a}{\gamma }\right)
			+\arctan\left(\frac{2 a}{\sigma ^2}\right)}{a \left(2 \gamma +\sigma ^2\right)}. \]
\end{example}
\begin{example}
	\label{ex:counter1}
	\rm{Consider  the example
		$g(\d \omega) = \frac{1}{2} \delta_{\omega_0}(\d \omega) + \frac{1}{2} \delta_{-\omega_0}(\d \omega)$
		for $\omega_0 > 0$. Set  $\beta = 1$, $\sigma = 1$, and $\omega_0 = 2$, then
		with these parameters, $\kappa_c(g) \approx 11.18$.  However,
		for $\kappa = 9$, Figure~\ref{fig:Gkappa} shows that
		$G_\kappa$ has three fixed points.  Hence, in this example
		the upper bound given by
		$\kappa_c(g)$ is not sharp.  Example \ref{ex:counter2} below
		further studies this example.
		\begin{figure}[ht]
			\centering
			\vspace{-1em}
			\input{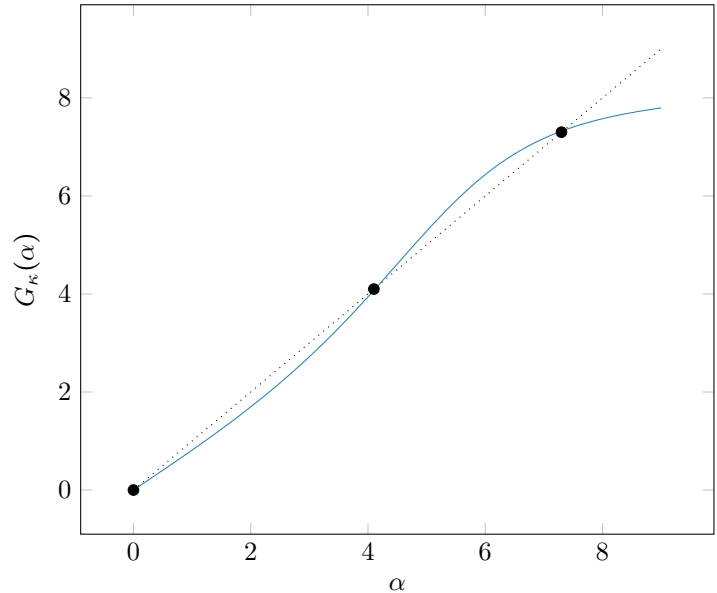}
			\caption{The function $G_\kappa$ for $\kappa = 9$, $g(\d \omega) = \frac{1}{2} \delta_{\omega_0}(\d \omega) + \frac{1}{2} \delta_{-\omega_0}(\d \omega)$ with $\beta = 1$, $\sigma = 1$, and $\omega_0 = 2$. This function $G_\kappa$ has three fixed points, despite that $\kappa < \kappa_c(g) \approx 11.18$.
			}
			\label{fig:Gkappa}
		\end{figure}}
\end{example}

\vspace{-1em}
%%%%%%%%%%%%%%%%%%%%%%%%%%%%%%%%%%%%%%%%%%
\section{Stability of the uniform measure}
\label{s.stable}
\def\dist{\mathfrak{g}_m} % please change if you do not like.  To me $d_m$ was confusing
In this section, we limit ourselves to intrinsic frequency distributions $g$ that are symmetric around zero, and study the stability of the $g$-uniform Nash equilibrium, as a function of $g$, $\kappa$, $\beta$, and $\sigma$.

We start with a definition of (local) stability.
For any $\mu \in \cP(\T)$,  as in \cite{CCS}, we set:
\begin{equation}
	\label{eq:def-dist}
	\dist(\mu) := \max \left( \left| \int_\T \cos(x) \mu(\d x) \right|, ~  \left| \int_\T \cos(2x) \mu(\d x) \right|, \left| \int_\T \sin(x) \mu(\d x) \right|, \left| \int_\T \sin(2x) \mu(\d x) \right| \right).
\end{equation}
Let $m$ and $\unig$ be as in \eqref{eq:uniform}.
As $\dist(\unig)=0$, one may view $\dist(\nu)$ as a pseudo distance of $\nu$
to the $g$-uniform measure $\unig$, explaining the subscript $m$.
\begin{definition}
	\label{def:stable}
	We say that the $g$-uniform measure $\unig$ is {\rm{stable}} at $\kappa$,
	if there exist constants $C,\lambda, \epsilon >0$ such that for any
	initial condition $\nu(\d x, \d \omega) = \nu^\omega(\d x) g(\d \omega)$
	satisfying
	$$
		\dist(\nu^\omega) \leq \epsilon, \qquad g(\d \omega)-a.s.
	$$
	there exists a solution
	$(\mu_t)_{t\ge 0}$ of the MFG problem with parameter $\kappa$
	with $\mu_0=\nu$, satisfying
	$$
		\dist(\mu^\omega_t) \leq C e^{-\lambda t} \dist(\nu^\omega),
		\qquad   \forall t \geq 0, \quad g(\d \omega)-a.s.,
	$$
	where $\mu^\omega$ is given by
	$\mu_t(\d x, \d \omega) = \mu^\omega_t(\d x) g(\d \omega)$.
\end{definition}
\begin{remark}
	If stability holds with $\dist$ as above, then one can control the distance to $m$ in stronger norms, for instance any negative Sobolev norms: for any $s < 1/2$, there is a constant $C_s$ such that:
	\[ \norm{\mu^\omega_t-m}_{H_{s}(\T)} \leq C_s e^{-\lambda t }
		\norm{\nu^\omega - m}_{H_{s}(\T)}, \quad \forall t \geq 0,
		\quad g(\d \omega) - a.s.   \]
	In fact, $\dist$ measures all the information relevant for assessing stability.
\end{remark}
In the subsequent subsections, we prove Theorem~\ref{th:lb1}
and Theorem~\ref{th:lb2}.

\subsection{Reformulation}
\label{ss.reformulation}
For $\lambda >0$,  set
$$
	L^\infty_\lambda:= \{ k: \R_+ \rightarrow \R \ :\ \norm{k}^\infty_\lambda<\infty\},
	\quad \text{where} \quad
	\norm{k}^\infty_\lambda := \text{esssup}_{t \geq 0} \abs{k(t)} e^{\lambda t}.
$$
%\rc{We may want to define here $\|r\|^1_\lambda$ which we use later on. We also use the norm $\|\phi\|^\infty$ and I do not know if we define it anywhere.}
Given $h = (h^1, h^2) \in L^\infty_\lambda \times L^\infty_\lambda$, we consider the cost
$$
	\ell^h(t, x) = -h^1(t) \cos(x) - h^2(t) \sin(x),
$$
and consider the following optimization problem
for a player with intrinsic frequency
\[ v^\omega(t, x) = \inf_{\alpha \in \cA_t}  \E \int_t^\infty e^{-\beta(u-t)}
	\left( \frac{1}{2} \alpha^2_u + \ell^h(u, X^{\omega, \alpha}_u) \right) \d u ,  \]
where $X^{\omega, \alpha}_u = x + \omega (u-t) + \int_t^u \alpha_\theta \d \theta + \sigma (B_u - B_t)$.
The corresponding value function is the unique smooth solution of the HJB equation
$$
	\partial_t v^\omega  - \frac{1}{2} (\partial_x v^\omega)^2
	+ \frac{\sigma^2}{2} \partial_{xx} v^\omega + \ell^h + \omega \partial_x v^\omega = \beta v^\omega.
$$
We now fix an initial condition $\mu_0=\nu=\nu^\omega(\d x) g(\d \omega)$,
and let $X^{\omega, *}_t$ be the state of the player following the optimal strategy. Then,
$$
	X^{\omega,*}_t = X^\omega_0 + \omega t  - \int_0^t \partial_x v^\omega_s(X^{\omega,*}_s) \d s + \sigma B_t,
$$
with initial condition $\cL(X^\omega_0) = \nu^\omega$. Set
$$
	\Phi(\nu, h) := - h + \kappa
	\begin{pmatrix}  \int \E \cos(X^{\omega,*}_t) g(\d \omega) \\ \int{\E \sin(X^{\omega,*}_t) g(\d \omega)}  \end{pmatrix}_{t \geq 0},
	\quad \nu \in \cP(\T \times \R), \ h \in   L^\infty_\lambda \times L^\infty_\lambda.
$$
{
We equip the space $\cP(\T \times \R)$ with the following gauge function.
For a given probability measure $\nu$ on $\mathbb{T} \times \mathbb{R}$ with second marginal $g$, the Polish structure of the space ensures by the Disintegration Theorem the existence of a $g$-almost everywhere unique family of conditional probabilities $\{\nu^\omega\}_{\omega \in \mathbb{R}}$ on $\mathbb{T}$ such that $\nu(\d x, \d \omega) = \nu^\omega(\d x)g(\d \omega)$. For any probability measure $\mu$ on $\mathbb{T}$, recall that $\dist(\mu)$ is defined in~\eqref{eq:def-dist}.
The function $\omega \mapsto \dist(\nu^\omega)$ is well-defined and Borel measurable. Set
\[ \delta(\nu) = g\text{-ess sup}_{\omega \in \mathbb{R}} \, \dist(\nu^\omega). \]
Because the disintegration is unique up to $g$-null sets, this quantity is well-defined and independent of the choice of conditional probabilities.
}
The above map $\Phi$ allows us to
construct solutions to the MFG problem.
Indeed, it is clear that  if $\Phi(\nu, h) = 0$
for some $h \in L^\infty_\lambda$, then
$$
	\mu_t(\d x, \d \omega) := \cL(X^{\omega, *}_t) g(\d \omega)
$$
is a solution of the MFG problem satisfying $\mu_0 = \nu$.
Moreover, since a function $h(\nu)\in L^\infty_\lambda \times L^\infty_\lambda$ satisfying $\Phi\bigl(\nu,h(\nu)\bigr)=0$ decays towards $0$ at an exponential rate, we then deduce that the corresponding MFG solution converges to $\unig$ as $t$ goes to infinity.
Therefore, if there are solutions $\Phi(\nu, h) = 0$
for any $\nu$ close enough to $\unig$, then we
can conclude that the $g$-uniform distribution is stable
in the sense of Definition \ref{def:stable}.
To construct such solutions, we observe that the
$g$-uniform distribution satisfies $\Phi(\unig,0)=0$
and use the implicit function theorem.

We start our analysis by computing the derivative of $\Phi$ at the uniform solution.
The following is proved in Section~\ref{sec:Frechet} below.

\begin{lemma}
	\label{lem: Frechet}
	Suppose that $g$ is symmetric around the origin. Then,
	$\Phi$ is Fr\'echet differentiable at $(\unig,0)$ and
	for any $h = (h^1, h^2) \in L^\infty_\lambda \times L^\infty_\lambda$:
	\[
		D_h \Phi(\unig, 0)\cdot h
		= -h + \kappa \begin{pmatrix} L h^1 \\ L h^2 \end{pmatrix},
	\]
	where for any $k \in L^\infty_\lambda$
	\begin{equation}
		\label{eq:operator-L}
		(L k)(t)= \frac{1}{2} \int_0^t \int_\theta^\infty e^{-\frac{\sigma^2}{2}(t-\theta)}
		e^{-\gamma(u-\theta)} \int g(\d \omega) k(u) \cos(\omega(t-u))\ \d u \ \d \theta,
		\qquad t \ge 0.
	\end{equation}
\end{lemma}

Next, we obtain a sufficient condition for the solvability of the equation $\Phi(\nu, h) = 0$ in $h \in L^\infty_\lambda \times L^\infty_\lambda$.
Let  $I$ be the identity operator.

\begin{proposition}
	\label{pro:sufficient}
	Suppose that the linear operator
	$I-\kappa L$ is invertible on $L^\infty_\lambda$
	for some $\lambda>0$.  Then, the
	$g$-uniform distribution $\unig$ is stable at $\kappa$.
	In particular,
	$$
		\kappa_2(g) \ge \frac{1} {\norm{L}^\infty_\lambda},
		\qquad \text{where} \qquad
		\norm{L}^\infty_\lambda := \sup \{ \norm{ L k}^\infty_\lambda :  \norm{k}^\infty_\lambda = 1\}.
	$$
\end{proposition}
\begin{proof}
	The analysis of Section~\ref{sec:estimates}
	implies that  $\Phi$ and its derivatives are
	continuous.  Then, when  $I-\kappa L$ is invertible,
	Lemma~\ref{lem:into} and the continuity results  imply that the
	Implicit Function Theorem~\ref{th:IFT} applies to the equation $\Phi(\nu,h)=0$.
	Thus, there is $\eta>0$ such that
	for any $\nu = \nu^\omega(\d x) g(\d \omega)$ with $\dist(\nu^\omega) \le \eta$,
	there exists $h (\nu)\in L^\infty_\lambda$ such that
	$\Phi(\nu,h(\nu))=0$.  As argued before, this solution $h(\nu)$
	is related to a Nash equilibrium which converges exponentially
	to the $g$-uniform.  Hence,  $\unig$ is stable at $\kappa$.

	To prove the final inequality, we observe that
	$I-\kappa L$ is invertible when
	$\kappa < 1/ \norm{L}^\infty_\lambda$.
\end{proof}

\subsection{Proof of Theorem~\ref{th:lb1}}
\label{ss.prooflb1}

Since $g$ is symmetric around zero,
$$
	\int g(\d \omega) \cos (\omega(t-u))  = \hat{g}(i(t-u)).
$$
Recall $\gamma=\sigma^2/2 + \beta$ of \eqref{eq:def-gamma} and
for $u,t \in \R_+$ set
$$
	K(t, u) := \int_0^{u \wedge t} e^{-\frac{\sigma^2}{2}(t-\theta)} e^{-\gamma(u-\theta)} \d \theta
	= \frac{ e^{-\frac{\sigma^2}{2}t} e^{-(\frac{\sigma^2}{2} + \beta) u}}{\sigma^2 + \beta} \left( e^{(\sigma^2+\beta) (u \wedge t)  } - 1 \right)
	=\frac{ e^{-\frac{\sigma^2}{2}t} e^{-\gamma u}}{\frac{\sigma^2}{2} + \gamma} \left( e^{(\frac{\sigma^2}{2}+\gamma) (u \wedge t)  } - 1 \right).
$$
By Fubini, the operator $L$ defined in \eqref{eq:operator-L} satisfies:
$$
	(L k)(t) = \frac{1}{2} \int_0^\infty K(t, u) \hat{g}(i(t-u)) k(u) \d u.
$$
Theorem~\ref{th:lb1} follows
from the following result and Proposition~\ref{pro:sufficient}.

\begin{lemma}
	\label{lem:first}
	Assume that $g$ is symmetric around zero. Then
	\[ \norm{L}^\infty_\lambda \leq \frac{ 1 } {(\gamma + \lambda)(\sigma^2 - 2 \lambda)}. \]
	In particular, if $\kappa < \kappa_c(\delta_0)=\gamma \sigma^2$ then there exists
	$\lambda > 0$ small enough such that $I - \kappa L$ is invertible as a linear operator
	in $\cL(L^\infty_\lambda; L^\infty_\lambda)$.
\end{lemma}
\begin{proof}
	Note that as $g$ is a probability measure, we have $\abs{\hat{g}(i t)} \leq 1$.
	The result then follows from the explicit expression of $K(t, u)$.
\end{proof}

\subsection{Proof of Theorem~\ref{th:lb2}}
\label{ss.prooflb2}

Recall $\kappa_c(g)$ of \eqref{fo:tilde_kappa_c}.

\begin{lemma}
	\label{lem:second}
	Assume that $\hat{g}(i t) \geq 0$ for all $t \geq 0$. It holds that for $\lambda \in (0, \sigma^2/2)$,
	\[ \norm{L}^\infty_\lambda =
		\int_{\R} \frac{ (\sigma^2 - 2 \lambda)(\lambda + \gamma) + 2 \omega^2}{\left[(\sigma^2 - 2\lambda)^2 + 4 \omega^2\right] \left[(\lambda + \gamma)^2 + \omega^2\right]}
		g(\d \omega). \]
	In particular, if $\kappa < \kappa_c(g)$, then there exists $\lambda > 0$ small enough such that $I - \kappa L$
	is invertible as a linear operator in $\cL(L^\infty_\lambda; L^\infty_\lambda)$,
	and consequently Theorem~\ref{th:lb2} holds.
\end{lemma}
\begin{proof}
	As $\hat{g}(i(t-u)) \geq 0$, we have
	\begin{align*}
		\norm{L}^\infty_\lambda & = \sup_{t \geq 0} \frac{1}{2}
		\frac{e^{-\frac{\sigma^2}{2} t} e^{\lambda t}}{\sigma^2 + \beta} \int_0^\infty e^{-(\lambda + \gamma)u} \left(e^{(\sigma^2 + \beta) u \wedge t} -1 \right) \hat{g}(i (t-u)) \d u
		\\
		                        & =  \sup_{t \geq 0} \frac{1}{2} \frac{e^{-\frac{\sigma^2}{2} t}
			                                                         e^{\lambda t}}{\sigma^2 + \beta} \int_{\R} g(\d \omega)
		\int_0^\infty e^{-i \omega (t-u)}e^{-(\lambda + \gamma)u} \left(e^{(\sigma^2 + \beta) u \wedge t} -1 \right) \d u.
	\end{align*}
	We then utilize the following identity,
	valid for all $A, B \in \C$ with $\Re(A) > 0$ and $B \neq A$:
	$$
		\int_0^\infty e^{-A u} (e^{B (t \wedge u)} - 1) \d u = \frac{B (e^{(B-A) t} -1)}{A(B-A)}.
	$$
	We choose $A = -i \omega + \lambda + \gamma$ and  $B = \sigma^2 + \beta$.
	Then,  $B-A = \sigma^2/2 + i \omega -\lambda$ and
	\begin{align*}
		\frac{1}{\sigma^2 + \beta}\ \int_0^\infty e^{-i \omega (t-u)}e^{-(\lambda + \gamma)u}
		\left(e^{(\sigma^2 + \beta) u \wedge t} -1 \right) \d u
		 & =e^{-i \omega t} \frac{ e^{(\sigma^2/2 + i \omega - \lambda)t} - 1 }
		                    { (\sigma^2/2 +i \omega - \lambda )(\lambda + \gamma - i \omega)}.
	\end{align*}
	We use this and the fact that $g$ is symmetric to arrive at
	\begin{align*}
		\norm{L}^\infty_\lambda & =\sup_{t \geq 0} \frac{1}{2}
		\int_{\R} g(\d \omega) e^{-i \omega t} \frac{ 1- e^{-(\sigma^2/2 + i \omega - \lambda)t}}
		                                       { (\sigma^2/2 +i \omega - \lambda )(\lambda + \gamma - i \omega)} \\
		                        & = \int_{\R} g(\d \omega)
		\frac{ (\sigma^2 - 2 \lambda)(\lambda + \gamma) + 2 \omega^2}
		{\left[(\sigma^2 - 2\lambda)^2 + 4 \omega^2\right]
			\left[(\lambda + \gamma)^2 + \omega^2\right]} - C:= \kappa(\lambda)- C,
	\end{align*}
	where
	$$
		C = \inf_{t \geq 0}  \frac{1}{2} e^{(\lambda -\frac{\sigma^2}{2}) t}
		\int_{\R}e^{-i \omega t} \varphi(i \omega) g(\d \omega),
		\quad \text{ and } \quad
		\varphi(i \omega) := \frac{1}{ (\sigma^2/2 +i \omega - \lambda )(\lambda + \gamma - i \omega)}.
	$$
	The integral, $ \int_{\R} g(\d \omega)
		e^{-i \omega t} \varphi(i \omega) g(\d \omega)$ is the Fourier transform of the product of $g$  and $\varphi$.\
	Moreover, by a direct computation,
	we conclude that the Fourier transform of $\varphi(i \omega)$ is non-negative and  the Fourier transform of $g$ is non-negative
	by hypothesis, the convolution of two non-negative
	functions is also non-negative.  As $\lambda < \frac{\sigma^2}{2}$, by
	letting $t$ tend to infinity in the above minimization, we conclude that $C = 0$.
	Additionally, using that $2\gamma= 2\beta +\sigma^2$, we directly
	calculate that $\kappa'(0)<0$.
	This implies that for all $\kappa \le \kappa_c(g)$,
	$\kappa \norm{L}^\infty_\lambda <1$
	provided that $\lambda$ is sufficiently small.
	Then, Proposition~\ref{pro:sufficient}
	implies that $\unig$ is stable for these $\kappa$
	values completing the proof of Theorem~\ref{th:lb2}.
\end{proof}

\section{Penrose Characterization}
\label{sec:Penrose}

In view of Proposition~\ref{pro:sufficient}, the stability of $\unig$
is deduced from the invertibility of $I - \kappa L$.
For general distributions, Laplace transform techniques
can be used to obtain sharper lower  bounds for $\kappa_2(g)$.
In the next section, we leverage this observation
to study the distributions that are the sum of two point masses.

Recall, the strip $\mfS$ and the function $P$ defined in \eqref{eq:P}
and the Penrose constant $\kappa_P(g)$ of~\eqref{eq:ckappa-penrose}.
The following properties of $P$
derived from the symmetry of $g$ are used in the sequel.
\begin{lemma}
	\label{lem:symmetry}
	$P$ is holomorphic on $\mfS$, $P(\bar{z}) = \overline{P(z)}$,
	and it is symmetric with respect to $\beta/2$:
	$$
		P(z + \beta/2) = P(-z + \beta/2), \qquad \forall z \in \mfS.
	$$
\end{lemma}
\noindent
The location of the zeros of $P$ encodes information on the stability of
$\unig$ and we formulate this as follows.
\vspace{5pt}

\noindent
{\bf{Penrose condition.}}
{\emph{Assume that there exists $\lambda \in (0, \sigma^2/2)$ such that
		\begin{equation}
			\label{fo:sufficient}
			1 - \frac{\kappa}{2} P(z) \neq 0,\qquad
			\text{whenever} \ \ \Re(z) \in [-\lambda, \beta+\lambda].
		\end{equation}}}
\noindent
Similarly to the standard Kuramoto dynamical system~\cite{Dietert16}, there is a geometric way to check the condition.
Indeed, consider the curve $\{ P(i \theta), \ \theta \in \R\}$
in the complex plane.
As $\lim_{ \theta \rightarrow \pm \infty} P(i \theta) = 0$, this curve
starts and ends at $z = 0$.
One can show that for small values $\kappa < \gamma \sigma^2$, $2/\kappa$
is outside this curve.  This and the definition
of $\kappa_P(g)$ prove the following.
\begin{lemma}
	\label{lem:kappaP}
	The Penrose condition \eqref{fo:sufficient} holds
	for all  $\kappa< \kappa_P(g)$ for all $g$ with a  finite second moment.
\end{lemma}
\begin{proof}
	Using the estimates of Lemma~\ref{lem:estimate-decay-imz} below and the argument principle, we deduce that the number of zeros of $1- \frac{\kappa}{2} P(z)$ within the strip $\Re(z) \in [0, \beta]$ is equal to:
	\[
		G(\kappa) :=  \lim_{\theta \uparrow +\infty } -\frac{1}{2 \pi i } \int_{\Gamma_\theta} \frac{P'(z) }{1 - \frac{\kappa}{2} P(z)} \d z,
	\]
	where $\Gamma_\theta$ is the simple anticlockwise contour of the rectangle with corners $\{- i \theta, \beta - i \theta,  \beta + i \theta, + i \theta \}$.
	Indeed, for $\theta$ large enough, if $\abs{\Im(z)} \geq \theta$ then $\frac{\kappa}{2} \abs{P(z)} < 1$, and so $1 - \frac{\kappa}{2} P(z)$ is not equal to zero. So the number of zeros within the strip is the same as the number of zeros within $\Gamma_\theta$. In addition, the integrals along the horizontal sides vanish as $\theta \uparrow +\infty$. Using the symmetry of Lemma~\ref{lem:symmetry}, one deduces that
	\[
		G(\kappa) = \frac{1}{\pi} \int_{-\infty}^\infty \frac{P'(i y)}{1 - \frac{\kappa}{2} P(i y)} \d y.
	\]
	As $\kappa_P(g) = \inf\{ \kappa \geq 0: \exists y \in \R \text{ s.t. } 1-\frac{\kappa}{2} P(i y) = 0 \}$,
	it holds that $G$ is continuous on $[0, \kappa_P(g))$. Because $G$ is integer-valued,
	this implies that  $G$ is constant and is equal to $G(0) = 0$ for all $\kappa \in [0, \kappa_P(g))$.
	Since the zeros of $1- \frac{\kappa}{2} P(z)$ are isolated in $\mfS$
	and there is no zero within the strip $\Re(z)\in[0, \beta]$,
	we conclude that for all  $\lambda > 0$ small enough
	there are no zeros within the strip $[-\lambda, \beta + \lambda]$ as well.
	Therefore, the Penrose condition holds.
\end{proof}
\begin{lemma}
	\label{lem:estimate-decay-imz}
	There is a constant $C(\beta, \sigma)$ such that for all $z \in \C$:
	\[ \Re(z) \in [0, \beta] \implies \abs{P(z)} + \abs{P'(z)} \leq  \left( 1 + \int_\R \omega^2 g(\d \omega) \right) \frac{C(\beta, \sigma)}{\abs{\Im(z)}^2}.  \]
\end{lemma}
\begin{proof}
	Using the explicit expression for $P(z)$, we find that for $z = x + i y$,
	\begin{align*}
		P(z) & = \int_\R \frac{g(\d \omega)}{ (\gamma-x)(\frac{\sigma^2}{2}+x) + (\omega-y)^2 + i(\omega-y) (2x-\beta) }                                                                                                  \\
		     & = \int_\R \frac{ (\gamma-x)(\frac{\sigma^2}{2}+x) + (\omega-y)^2 - i(\omega-y) (2x-\beta)  }{  \left((\gamma-x)(\frac{\sigma^2}{2}+x) + (\omega-y)^2\right)^2  + (\omega-y)^2 (2x-\beta)^2 } g(\d \omega).
	\end{align*}
	Let $A(x) = (\gamma - x)(\frac{\sigma^2}{2} +x) \in [\gamma \frac{\sigma^2}{2}, \frac{(\sigma^2 +\beta)^2}{4}]$.
	We use the following inequality:
	\begin{equation}
		\label{eq:foupperbound}
		\frac{y^2}{A(x)+(\omega-y)^2} \leq 1 + \frac{\omega^2}{A(x)},
	\end{equation}
	to deduce that
	\[ \abs{\Re(P(z))} \leq \left( 1 + \frac{2}{\sigma^2 \gamma} \int_\R \omega^2 g (\d \omega) \right) \frac{1}{\abs{\Im(z)}^2}.   \]
	To treat the imaginary part, we first use that $2x-\beta \leq \beta$ and then rely on the following inequality:
	\[ \frac{\alpha}{A(x) + \alpha^2} \leq \frac{1}{\sqrt{2 A(x)}}, \quad \text{ where } \alpha = \abs{\omega - y}. \]
	We therefore find, using again \eqref{eq:foupperbound}:
	\[ \abs{\Im(P(z))} \leq \frac{\beta}{\sqrt{2 A(x)}} \left( 1 + \int_\R \frac{\omega^2}{A(x)} g(\d \omega) \right) \frac{1}{y^2}. \]
	As $A(x) \geq \gamma \frac{\sigma^2}{2}$, one deduces the result for $P(z)$.
	For the derivative:
	\begin{equation*}
		P'(z)=\int_{\R}\frac{ g(\d \omega)}{(\gamma+i\omega-z)(\sigma^2/2+z-i\omega)}\Bigl(
		\frac{1}{\gamma+i\omega-z}-\frac{1}{\sigma^2/2+z-i\omega}\Bigr).
	\end{equation*}
	One can now use the bound found for $P(z)$ because
	$$
		\Bigl|\frac{1}{\gamma+i\omega-z}\Bigr|\le \frac{3}{2(\gamma-x)} \quad \text{ and } \quad
		\Bigl|\frac{1}{\sigma^2/2+z-i\omega}\Bigr|\le \frac{3}{2(\sigma^2/2+x)}.
	$$
	This completes the proof.
\end{proof}

\begin{figure}
	\centering
	\input{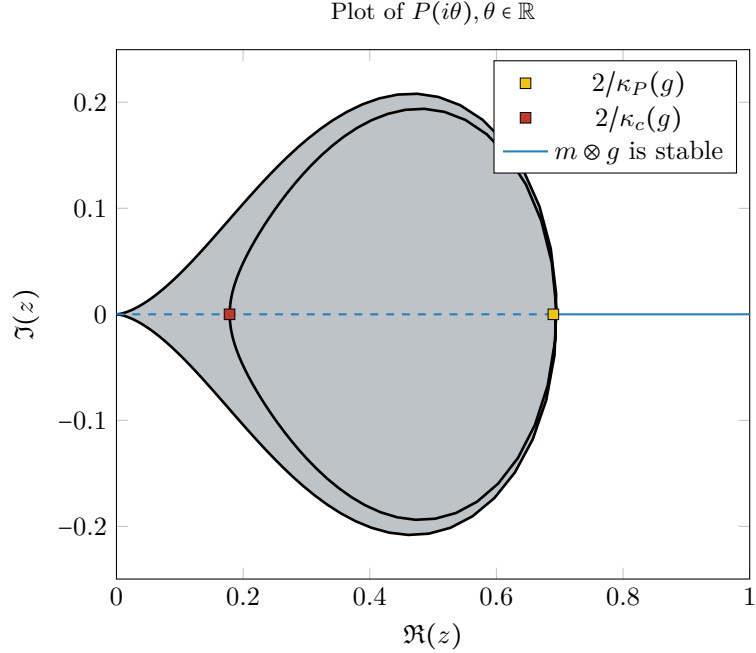}
	\caption{The Penrose diagram for two Dirac measures. The critical parameter $2/\kappa_P(g)$ is at the rightmost intersection between the curve $\{P(i \theta), \theta \in \R \}$ and $\Im(z) = 0$.
		Our result shows that if $\kappa$ is such that $2/\kappa$ belongs to the exterior of the Penrose curve, then $\unig$ is stable. In this example, the two critical parameters do not match: $\kappa_P(g) < \kappa_c(g)$.
	}
	\label{fig:penrose}
\end{figure}

When $g$ is the Dirac measure at zero,
$P(i\theta)= (\gamma \frac{\sigma^2}{2} +\theta^2 + i \beta \theta)^{-1}$
is real only for $\theta=0$ yielding
$\kappa_P =\kappa_c$
consistent with the results of~\cite{CCS}.
We believe that this identity  holds
under the hypothesis of Theorem~\ref{th:lb2}:
$$
	\mathfrak{F}g(t) \ge 0, \  \forall t \in \R
	\qquad \Rightarrow  \qquad
	\kappa_P(g) = \kappa_c(g).
$$
However, in general these two
thresholds do not match, as demonstrated in the
following example.
\begin{example}
	\label{ex:counter2}
	{\rm{Consider the case $\omega_0 = 2$, $g(\d \omega) = \frac{1}{2} \delta_{-\omega_0}(\d \omega) + \frac{1}{2} \delta_{\omega_0}(\d \omega)$, $\beta = 1$ and $\sigma = 1$.
			For this example, we have $\kappa_P(g) \approx 2.9 < \kappa_c(g) \approx 11.2$. The associated Penrose diagram is illustrated in Figure~\ref{fig:penrose}. }
	}
\end{example}
\subsection{Formulation in the Laplace domain}
\label{ss.Laplapce}

For $h \in L^\infty_\lambda$,
$$
	\hat{h}(z) = \int_0^\infty e^{-zt } h(t) \d t \in \R,
	\qquad \text{when} \ \Re(z) > -\lambda.
$$
\begin{lemma}
	Suppose $h \in L^\infty_\lambda$
	for some $\lambda \in (0, \sigma^2/2)$.
	Then
	$$
		\widehat{Lh}(z) = \frac{1}{2} \int_{\R} \frac{ \hat{h}(z)
			- \hat{h}(\gamma + i \omega)}{(z+\sigma^2/2 - i \omega)(\gamma + i \omega-z)}
		g(\d \omega),\qquad
		\qquad \text{for all} \ z \in \C
		\ \text{with}\ \ \Re(z) > -\lambda.
	$$
\end{lemma}
\begin{proof}
	By Fubini, we have
	$$
		\widehat{Lh}(z) = \frac{1}{2} \int_{\R} \int_0^\infty \int_0^\infty e^{-zt }K(t, u)
		e^{i \omega (t-u) } \d t k(u)\d u g(\d \omega).
	$$
	A direct calculation shows that
	$$\
		\int_0^\infty e^{-zt }K(t, u) e^{i \omega (t-u) } \d t
		= \frac{e^{-z u} - e^{-(\gamma + i \omega) u}}{(z+\sigma^2/2 - i \omega)(\gamma + i \omega -z)}.
	$$
	Hence, the result follows.
\end{proof}

We next obtain an equivalent formulation of the equation $(I - \kappa L)h = \phi$.
\begin{lemma}
	\label{lem:Laplace}
	Assume that \eqref{fo:sufficient} holds for
	some  $\lambda \in (0, \sigma^2/2)$.
	Suppose that $h, \phi \in L^\infty_\lambda$ satisfy the equation
	$(I - \kappa L)h = \phi$.
	Then, for all $z \in \C$ with $\Re(z) > -\lambda$ the following
	equation holds:
	\begin{equation}
		\label{eq:master-eq-laplace}
		\hat{h}(z) = \hat{\phi}(z) + \frac{\kappa}{ 2- \kappa P(z)} \int_{\R}
		\frac{\hat{\phi}(z) - \hat{h}(\gamma + i \omega)}{(z + \sigma^2/2 - i \omega)(\gamma + i \omega - z)}
		g(\d \omega). \end{equation}
\end{lemma}
\begin{proof}
	The equation $(I - \kappa L) h = \phi$ writes
	\[ \hat{h}(z)(1 - \frac{\kappa}{2} P(z)) + \frac{\kappa}{2} \int_{\R} \frac{\hat{h}(\gamma + i \omega)}{(z+\sigma^2/2 - i \omega)(\gamma + i \omega-z)} g(\d \omega)  = (1- \frac{\kappa}{2} P(z)) \hat{\phi}(z) + \frac{\kappa}{2} P(z) \hat{\phi}(z).  \]
	Dividing by $1 - \frac{\kappa}{2} P(z)$, we deduce the result.
\end{proof}
{
\section{The case of symmetric Dirac point masses}

\label{s.two_Diracs}
Let $\omega_1, \cdots \omega_N$ be distinct strictly positive real numbers. 
We set $\omega_{-i} = -\omega_i$ and let
\[ I = \{-N, -(N-1), \cdots, -1, 1, \cdots, N-1, N \}. \]
In this section, we consider  intrinsic frequency distributions given by
\[
	g(\d \omega) = \sum_{i \in I} \alpha_i \delta_{\omega_i}(\d \omega),
\]
where $\alpha_{-i} = \alpha_i > 0$ with $\sum_{i \in I}\alpha_i =1$.
The main result of this section is the following:
\begin{proposition}
	\label{prop:inversibility-I-L}
	Suppose that the Penrose condition
	\eqref{fo:sufficient} holds for some $\lambda >0$.
	Then, for any $\phi \in L^\infty_\lambda(\R_+)$
	there is a unique $k \in L^\infty_\lambda(\R_+)$ such that
	$$
		k(t) = \phi(t) +\kappa L k(t), \quad t \in \R_+.
	$$
	In particular,  the operator
	$I- \kappa L: L^\infty_\lambda \rightarrow L^\infty_\lambda$ is invertible with
	a bounded inverse.
\end{proposition}
Our main result, Theorem~\ref{th:two_Diracs} follows from this result:\\ \noindent{\emph{Proof of Theorem~\ref{th:two_Diracs}.}}
Fix $\kappa < \kappa_P(g)$.
In view of Lemma~\ref{lem:kappaP}
the Penrose condition holds, and by the
above lemma, $I-\kappa L$ is invertible.
We then invoke Proposition~\ref{pro:sufficient} and conclude
that the $g$-uniform distribution is stable.
Hence, $\kappa_2(g)$ is greater than $\kappa$.
\qed

We split the proof of Proposition~\ref{prop:inversibility-I-L} into two main parts: 
first invertibility in the Laplace domain, and then in the temporal domain.
\subsection{Invertibility in Laplace domain}
Set
\begin{align}
	Q(z, \omega) & := (C + i \omega -z)(C - i \omega + z) \quad \text{ where } \quad C := \frac{\beta + \sigma^2}{2}, \label{eq:Q}               \\
	\cN(z)         & := \prod_{i \in I} Q(z, \omega_i) - \frac{\kappa}{2} \sum_{i \in I} \alpha_i \prod_{j \neq i} Q(z, \omega_j). \label{eq:defN}
\end{align}
Let $P$ be as in~\eqref{eq:P}.  Then,
\[  P(z) = \sum_{i \in I} \frac{\alpha_i}{Q(z - \frac{\beta}{2}, \omega_i)}, \]
and so:
\[
	1 - \frac{\kappa}{2} P(z) = \frac{\cN(z-\frac{\beta}{2})}{ \prod_{i \in I} Q(z-\frac{\beta}{2}, \omega_i)}.
\]
Clearly, the sets of zeros of $1- \frac{\kappa}{2}P$ and $\cN$ are the same,
and since we assume that the Penrose condition~\eqref{fo:sufficient} holds, $1 - \frac{\kappa}{2}P$ has no roots on the strip 
$[-\lambda, \beta + \lambda]$ for some small but strictly positive $\lambda$. 
As $\cN(z)$ is a polynomial of degree $4N$, and that $\cN(z) = \cN(-z)$, we deduce that 
$\cN$ has exactly $2N$ zeros (counted with their multiplicity) that satisfy: $\Re(z) > 0$.
Let these distinct zeros be $\zeta_1, \cdots, \zeta_K$ with respective multiplicities $m_1, \cdots, m_K$, such that $\sum_{k=1}^K m_k = 2N$. Note that this set is closed under complex conjugation and may include purely real positive numbers.

Let $M$ be the $2N \times 2N$ matrix whose columns are indexed by $j \in I$, and whose rows are grouped by the roots $\zeta_k$. For each root $\zeta_k$, there are $m_k$ rows corresponding to the derivatives of order $r = 0, 1, \cdots, m_k - 1$ evaluated at $\zeta_k$:
$$M_{(k, r), j} = \frac{\partial^r}{\partial z^r} \left[ \frac{1}{Q(z, \omega_j)} \right]_{z = \zeta_k}$$
For example, if $\zeta_1$ has multiplicity 2, its first row is $\frac{1}{Q(\zeta_1, \omega_j)}$ and its second row is $- \frac{Q_z(\zeta_1, \omega_j)}{Q^2(\zeta_1, \omega_j)}$. If all the roots are simple, the matrix is simply given by
\[ M_{ij} = \frac{1}{Q(\zeta_i, \omega_j)}. \]
We first have:
\begin{lemma} The matrix $M$ is invertible.
\end{lemma}
\begin{proof}
	\textbf{Step 1}.  Assume that $M$ is not invertible. Then there exists a non-zero vector $x = (x_j)_{j \in I} \in \mathbb{C}^{2N}$ such that $Mx = 0$.
	Define the rational function $H(z) = \sum_{j \in I} x_j \frac{1}{Q(z, \omega_j)}$. The condition $Mx = 0$ is equivalent that for each root $\zeta_k$ of $\cN(z)$ in the right half-plane, $H(z)$ satisfies:
	\[ H^{(r)}(\zeta_k) = 0 \quad \text{for } r = 0, 1, \cdots, m_k - 1. \]
	This means $H(z)$ has a zero at $\zeta_k$ of multiplicity at least $m_k$.

	We note the conjugate symmetry of $Q(z, \omega)$: $\overline{Q(\bar{z}, \omega_j)} = Q(z, -\omega_j) = Q(z, \omega_{-j})$.
	Because the roots $\zeta_k$ and their multiplicities are symmetric with respect to the real axis (since $\cN(z)$ has real coefficients), the entire system of equations $Mx=0$ is closed under conjugation. Specifically, conjugating the system shows that the vector $y$ defined by $y_m = \overline{x_{-m}}$ is also in the null space of $M$.  Thus, without loss of generality, we assume there is a non-zero vector $c \in \ker M$ such that
	\[ c_j = \overline{c_{-j}}, \quad  \forall  j \in I. \]
	This symmetry guarantees that our newly formed rational function $H(z) = \sum_{j \in I} c_j \frac{1}{Q(z, \omega_j)}$ has real coefficients, i.e.~$\overline{H(\bar{z})} = H(z)$.

	\textbf{Step 2.} Define the function
	\[ W(z) = 1 -  \frac{\kappa}{2}\sum_{j \in I} \frac{\alpha_j}{Q(z, \omega_j)}. \]
	Using the common denominator $D(z) = \prod_{k \in I} Q(z, \omega_k)$, we can rewrite $W(z) = \frac{\cN(z)}{D(z)}$.
	For purely imaginary inputs $z = i\omega$, $Q(i\omega, \omega_j) = C^2 + (\omega - \omega_j)^2 > 0$. Because $N(z)$ has no zeros on the imaginary axis and $W(+\infty) = 1$, it follows that $W(i\omega) > 0$ for all $\omega \in \mathbb{R}$.
	Consider now:
	\[ F(z) = \frac{H(z)}{W(z)}. \]
	$F(z)$ can be expressed as $F(z) = H(z) \frac{D(z)}{\cN(z)}$. The poles of $F(z)$ in the right half-plane can only occur at the zeros of $\cN(z)$. By hypothesis, the roots of $\cN(z)$ in the right half-plane are exactly the points $\zeta_k$ with multiplicity $m_k$. However, since $c \in \ker M$, the numerator $H(z)$ vanishes at every $\zeta_k$ with multiplicity at least $m_k$.
	Consequently, the zeros of $H(z)$ cancel the poles of $1/N(z)$ in the right half-plane, meaning $F(z)$ is analytic in the right half-plane $\Re z \ge 0$.

	We now evaluate $F(z)$ at $z = C + i\omega_k$. At this point, $Q(C + i\omega_k, \omega_k) = 0$, so $D(z)$ vanishes. We therefore find:
	\[ F(C + i\omega_k) = - \frac{2}{\kappa}\frac{c_k}{\alpha_k}. \]
	\textbf{Step 3}. Consider the following positive definite integral along the imaginary axis:
	\[ I = \int_{-\infty}^\infty \frac{|H(i\omega)|^2}{W(i\omega)} d\omega \]
	Since $W(i\omega) > 0$, $I \ge 0$.
	Because $H(z)$ has real coefficients, $\overline{H(i\omega)} = H(-i\omega)$. Therefore, we can rewrite the integral as:
	$$I = \int_{-\infty}^\infty \frac{H(i\omega)}{W(i\omega)} H(-i\omega) \d \omega = \int_{-\infty}^\infty F(i\omega) H(-i\omega) \d \omega.$$
	Let $G(z) = F(z) H(-z)$. We evaluate $\int_{-\infty}^\infty G(i\omega) \d \omega$ by extending the contour into the right half-plane (RHP). The contour $\gamma$ goes up the imaginary axis and closes with a semicircle in the RHP. This orientation is clockwise so:
	\[ i \int_{-\infty}^\infty G(i\omega) \d \omega = -2\pi i \sum \text{Res}_{\text{RHP}} G(z). \]
	Because $F(z)$ is strictly analytic in the RHP, the only poles of $G(z)$ in the RHP come from $H(-z)$:
	\[ H(-z) = \sum_{j \in I} \frac{c_j}{C^2 - (z + i\omega_j)^2} \]
	The poles in the RHP are located where $(z + i\omega_j)^2 = C^2$, which yields $z = C - i\omega_j$.
	The residue of $H(-z)$ at $z = C - i\omega_j$ is evaluated as:
	\[ \lim_{z \to C - i\omega_j} (z - (C - i\omega_j)) \frac{c_j}{(C - z - i\omega_j)(C + z + i\omega_j)} = -\frac{c_j}{2C} \]
	At this pole, the value of the analytic part $F(z)$ is $F(C - i\omega_j) = F(C + i\omega_{-j}) = -\frac{2}{\kappa}\frac{c_{-j}}{\alpha_j} = -\frac{2}{\kappa} \frac{\overline{c_j}}{\alpha_j}$.
	Thus, the residue of $G(z)$ at $C - i\omega_j$ is:
	\[ \text{Res} = (-\frac{2 \overline{c_j}}{ \kappa\alpha_j}) \left(-\frac{c_j}{2C}\right) = \frac{|c_j|^2}{C \kappa \alpha_j}. \]
	Summing the residues over all $j \in I$, we evaluate the integral:
	\[ i \int_{-\infty}^\infty G(i\omega) \d \omega = -2\pi i \sum_{j \in I} \frac{|c_j|^2}{C \kappa \alpha_j} \implies I = -\frac{2\pi}{\kappa C} \sum_{j \in I} \frac{ |c_j|^2}{\alpha_j}.\]
	As $I \geq 0$, we deduce that $c_j = 0$ for all $j$; this contradicts the assumption that there is a non-trivial vector $c$ in the null space of $M$. Therefore,  $M$ is invertible.
\end{proof}
Next, we show that:
\begin{proposition}
	\label{prop:inversibility-Laplace-I-L}
	Let $\hat{\phi}$ be a holomorphic function on $\Re(z) > -\lambda$ satisfying
	\[ \overline{\hat{\phi}(z)} = \hat{\phi}(\bar{z}), \quad \Re(z) > -\lambda. \]
	There exists $2N$ constants $C(\gamma + i \omega_j) \in \C, \; j \in I$ such that:
	\begin{enumerate}
		\item	The function
		      \[ H(z) := \hat{\phi}(z) + \frac{\kappa/2}{ 1- \frac{\kappa}{2} P(z)} \int_{\R} \frac{\hat{\phi}(z) - C(\gamma + i \omega)}{Q(z - \frac{\beta}{2}, \omega)} g(\d \omega) \]
		      is holomorphic for $\Re(z) > -\lambda$.
		\item For all $j \in I$, $C(\gamma - i \omega_j) = \overline{C(\gamma + i \omega_j)}$.
		\item For all $j \in I$, $H(\gamma + i \omega_j) = C(\gamma + i \omega_j)$.
	\end{enumerate}
	Additionally, $\{C(\gamma + i \omega_j)\}_{\{j \in I\}}$ 
    are the unique constants which make $H$ holomorphic for $\Re(z) > -\lambda$.
\end{proposition}
\begin{proof}
	Note that
	\begin{align*}
		 & \int_{\R} \frac{\hat{\phi}(z) - C(\gamma + i \omega)}{(z + \sigma^2/2 - i \omega)(\gamma + i \omega - z)} g(\d \omega) = \sum_{j \in I} \alpha_j \frac{ \hat{\phi}(z) - C(\gamma + i \omega_j)}{ Q(z - \frac{\beta}{2}, \omega_j) } \\
		 & \qquad = \frac{ \sum_{ j \in I} \alpha_j (\hat{\phi}(z) - C(\gamma + i \omega_j)) \prod_{k \neq j} Q(z - \frac{\beta}{2}, \omega_k)}{ \prod_{k \in I} Q(z -  \frac{\beta}{2}, \omega_k)},
	\end{align*}
	so the function $H(z)$ can be rewritten as:
	\begin{equation}
		\label{eq:N-dirac-laplace}
		H(z) = \hat{\phi}(z) + \frac{\kappa}{2} \frac{  \sum_{ j \in I} \alpha_j (\hat{\phi}(z) - C(\gamma + i \omega_j)) \prod_{k \neq j} Q(z -  \frac{\beta}{2}, \omega_k) }{ \cN(z - \frac{\beta}{2})}.
	\end{equation}
	Our goal is to find the constants $\beta_j := C(\gamma + i \omega_j), j \in I$ which make $H(z)$ holomorphic for $\Re(z) > - \lambda$. Given $\beta \in \C^{2N}$, consider
	\[ F_\beta(z) = \sum_{ j \in I} \alpha_j (\hat{\phi}(z + \frac{\beta}{2}) - \beta_j) \prod_{k \neq j} Q(z, \omega_k), \quad \Re(z) > -\lambda. \]
	\textbf{Step 1}. We prove that there is a unique choice of $\beta \in (\C)^{2N}$ such that:
	\begin{equation} \label{eq:conditions-F-beta}  \forall k \in 1 \cdots K, \quad \zeta_k \text{ is a zero of order $m_k$ of $F_\beta(z)$}. \end{equation}
	Indeed, $F_\beta(z) = 0$ is equivalent to
	\[ \sum_{j \in I} \frac{\alpha_j \beta_j}{Q(z, \omega_j)} = \sum_{j \in I} \alpha_j \frac{\hat{\phi}(z + \beta/2)}{Q(z, \omega_j)}. \]
	Therefore, if for each $k \in 1\cdots K$, $\zeta_k$ is a zero of order $m_k$ of $F_\beta(z)$ then $\beta$ solves a linear system of the form
	\[ M \beta = x, \]
	for some vector $x \in \C^{2N}$. As $\det(M) \neq 0$, this linear system is uniquely solvable.
	We set:
	\[ \beta^*_j = \overline{\beta_{-j}}, \quad j \in I. \]
	\noindent \textbf{Step 2}. We now prove that $\beta^* = \beta$. First, note that $\overline{Q(z, \omega)} = Q(\bar{z}, -\omega)$. Therefore:
	\begin{align*}
		\overline{F_\beta(z)} & = \sum_{j \in I} \alpha_j (\hat{\phi}(\bar{z} + \frac{\beta}{2}) - \overline{\beta_j} ) \prod_{k \neq j} Q(\bar{z}, -\omega_k), \\
		                      & = \sum_{j \in I} \alpha_{-j} (\phi(\bar{z} + \frac{\beta}{2}) - \overline{\beta_j} )  \prod_{k \neq -j} Q(\bar{z}, \omega_k),   \\
		                      & = F_{\beta^*}(\bar{z}).
	\end{align*}
	Similarly, for all $p \geq 1$,
	\[ \overline{F^{(p)}_\beta}(z) = F^{(p)}_{\beta^*}(\bar{z}). \]
	So, if $\zeta_k$ is a zero of order $m_k$ of $F_\beta$, then $\overline{z_k}$ is a zero of order $m_k$ of $F_{\beta^*}$.
	As the set of zeros of $\cN(z)$ (and their order) is invariant by complex conjugation, we deduce that $\beta^*$ also satisfies \eqref{eq:conditions-F-beta}.
	By uniqueness, we deduce that $\beta^* = \beta$, which shows the second point of the proposition.

	\noindent \textbf{Step 3}. It remains to prove that $H(\gamma + i \omega_j) = C(\gamma + i \omega_j) = \beta_j$. To proceed, note that $Q(C + i \omega_k, \omega_k) = 0$. Therefore, \[ \cN(C + i \omega_k) = - \frac{\kappa}{2} \alpha_k \prod_{j \neq k} Q(C + i \omega_k, \omega_j). \]
	So, equation~\eqref{eq:N-dirac-laplace} gives (using that $\gamma - \beta/2 = C)$
	\begin{align*} H(\gamma + i \omega_j) & = \hat{\phi}(\gamma + i \omega_j) + \frac{\kappa}{2} \frac{  \alpha_j (\hat{\phi}(\gamma + i \omega_j) - \beta_j) \prod_{k \neq j} Q(C + i \omega_k - \frac{\beta}{2},\omega_k)  }{ - \frac{\kappa}{2} \alpha_j  \prod_{k \neq j} Q(C + i \omega_k - \frac{\beta}{2},\omega_k)  } \\
	                                      & = \beta_j.
	\end{align*}
	This completes the proof.
\end{proof}

\subsection{Invertibility in the temporal domain}
\label{ss.temporal}

In this subsection, we prove Proposition~\ref{prop:inversibility-I-L}. 
An outline of this technical proof given 
in Remark~\ref{rem:outline} below, clearly describes the various steps at a high level.
We invite the reader to review this outline before reading the proof in detail.

Recall that the operator $L$ is defined by
$$
	(Lk)(t) :=  \indica{t \geq 0 }\  \frac{1}{2} \int_0^\infty
	\left[ \int_{0}^{u \wedge t} e^{-\frac{\sigma^2}{2} (t-\theta)}
		e^{-\gamma (u-\theta)} \d \theta \right] \hat{g}(i(t-u)) k(u) \d u, \quad t \in \R.
$$
We have constructed a solution $\hat{h}$
of the Laplace transform of the equation $\phi =(I-\kappa L)h$.
We next show that this solution can be used to
obtain a solution in the temporal domain.

%Given $\lambda \in \R$ and $I$ a closed interval of $\R$, recall that $L^\infty_\lambda(I)$ denotes the space of measurable functions such that $\sup_{t \in I}$ $\abs{x(t)}e^{\lambda t} < \infty$. Similarly, $L^1_\lambda(I)$ denotes the space of functions such that $\int_I \abs{x(t)} e^{\lambda t} \d t < \infty$.
%\rc{We do what you want but I suggest that we remove the two preceding sentences. The notation $I$ conflicts with the identity operator and after all, I believe that we only use these function spaces over $I=\RR_+$ and $I=\RR$, and in those cases the notation is clear and self explanatory. Just a thought. }

For $k \in L_\lambda^\infty(\R)$, set
$$
	A k(t)  :=  \frac{\kappa}{2} \int_\R \left[ \int_{-\infty}^{u \wedge t}
		e^{-\frac{\sigma^2}{2} (t-\theta)} e^{-\gamma (u-\theta)}
		\ \hat{g}(i(t-u))\ \d \theta\right]\  k(u) \d u, \qquad t \in \R.
$$
By changing variables to $v = t-u$, we obtain
$$
	Ak(t) = \kappa \int_{\R} e^{-\frac{\sigma^2}{2} t}\  e^{-\gamma (t-v)}
	\ \frac{e^{(\gamma + \frac{\sigma^2}{2})(t \wedge (t-v))}}{ \sigma^2 + 2 \gamma} \hat{g}(i v) k(t-v)
	\ \d v.
$$
Next, we split the integral into $v \in \R_{-}$ and $v \in \R_{+}$,
$$
	Ak (t) = \int_{\R} a(v) k(t-v) \d v,
	\quad \text{where} \quad
	a(t) := \frac{\kappa}{\sigma^2+ 2 \gamma}\ \hat{g}(i t)\
	\left( e^{\gamma t} \indica{t < 0} + e^{-\frac{\sigma^2}{2} t }
	\indica{t \geq 0} \right), \quad t \in \R.
$$
Therefore, $A$ is a convolution operator and
$\hat{a}(z) = \frac{\kappa}{2}P(z)$ provided that $\Re(z) = -\lambda$.
Hence,
$$
	\widehat{A k}(z) = \frac{\kappa}{2} P(z) \hat{k}(z),
	\qquad \forall \ z \in \C \quad \text{with}\ \Re(z) = -\lambda.
$$

In addition, we have
\begin{lemma}
	\label{lem:invertibility-A}
	Let $\theta \in (-\gamma, \sigma^2/2)$. Suppose that
	$$
		1 - \frac{\kappa}{2} P( i \omega - \theta) \neq 0,
		\qquad
		\forall \omega \in \R.
	$$
	Then, for any $\phi \in L^\infty_\theta(\R)$
	there exists $k \in L^\infty_\theta(\R)$ such that
	$$
		k(t) = \phi(t) + (A k)(t), \qquad t \in \R.
	$$
\end{lemma}
\begin{proof}
	Set $a_\theta(t) := a(t) e^{\theta t} \in L^1(\R)$.
	We apply the Paley-Wiener theorem, \cite[Th. 4.3]{Grip}
	to construct $r_\theta \in L^1(\R)$ satisfying,
	$$
		r_\theta(t) = a_\theta(t) + \int_{\R} r_\theta(t-s) a_\theta(s) \d s
		= a_\theta(t) +(r_\theta * a_\theta)(t),
	$$
	where $*$ denotes the convolution.
	We write $r(t) = r_\theta(t) e^{-\theta t}$.
	{
		It holds that $r(t) \in L^1_{\theta}(\R)$, that is $t \mapsto r(t)
		e^{\theta t} \in L^1(\R)$ with
	\[ \norm{r}^1_\theta := \int_{\R} e^{\theta t} \abs{r(t)} \d t .
		\]
	}
In addition $r= a + a * r$.  Given $\phi \in L^\infty_\theta(\R)$, set
	$k := \phi+  r* \phi$.  Then,
	$$
		k= \phi+ r * \phi = \phi+ (a + a* r) * \phi =
		\phi+ a *(\phi + r * \phi )= \phi + a * k.
	$$
	Additionally, $\norm{k}^\infty_{\theta} \leq (1+\norm{r}^1_{\theta}) \norm{\phi}^\infty_{\theta} < \infty$.
	Hence,  $(I-A)^{-1} \phi = k$.

	To prove uniqueness observe that if $k=\phi+a*k$,
	then
	$$
		k=\phi+(r-a*r)*k=\phi +r*k-a*r*k=\phi+r*(k-a*k)=\phi+r*\phi.
	$$
	The open mapping theorem implies that
	this inverse is bounded.
\end{proof}

The connection between this new operator $A$ and the original operator
$L$ is given by the following.
\begin{lemma}
	Assume that $k \equiv 0$ on $\R_-$. Then for all $t \in \R$:
	\begin{align*} Ak(t) & = \kappa Lk(t) + \frac{\kappa \indica{t \geq 0}}{2}\int_{\R_+} \frac{e^{-\frac{\sigma^2}{2} t} e^{-\gamma s} }{\frac{\sigma^2}{2} + \gamma} \hat{g}(i(t-s)) k(s) \d s + \frac{ \kappa \indica{t < 0}}{2} \int_{\R_+} \frac{e^{\gamma (t-s)} }{\frac{\sigma^2}{2} + \gamma} \hat{g}(i(t-s)) k(s) \d s.
	\end{align*}
\end{lemma}
\begin{proof}
	This follows directly from the definition of $A$ and $L$, by distinguishing the cases where $t < 0$ and $t \geq 0$. When $t < 0$, by definition $Lk (t) = 0$, giving the result. When $t \geq 0$, as $u \geq 0$ as well, we can split the integral from $-\infty$ to $u \wedge t$ as an integral from $-\infty$ to $0$, plus an integral from $0$ to $u \wedge t$. This gives the stated decomposition.
\end{proof}

We now define for $t \in \R$,
$$
	(Bk)(t):=  - \kappa \indica{t \geq 0}\  \int_{\R_+}
	\frac{e^{-\frac{\sigma^2}{2} t} e^{-\gamma s} }{\sigma^2 +  2 \gamma} \
	\hat{g}(i(t-s)) k(s) \d s
	- \kappa \indica{t < 0}\ \int_{\R_+} \frac{e^{\gamma (t-s)} }{\sigma^2 +  2 \gamma}\
	\hat{g}(i(t-s)) k(s) \d s,
$$
such that
$$
	k \equiv 0 \ \text{ on }\ \R_-
	\quad \Rightarrow \quad
	\kappa L k  (t) = (A+B)k(t), \qquad t \in \R.
$$
Recall that $\hat{g}(it)$ is defined by \eqref{eq:def-fourier-transform-g}. Using this definition and the fact that $g$ is symmetric, we can rewrite $B k$ as follows:
\begin{align}
	\nonumber
	B k(t)
	 & = -\frac{\kappa }{\sigma^2 + 2 \gamma}
	\left[ e^{-\frac{\sigma^2}{2} t} \indica{t \geq 0}
		+  e^{\gamma t} \indica{t < 0} \right] \int_{\R} e^{i \omega t }\int_{\R_+}
	e^{-(\gamma + i \omega) s} k(s) \d s g(\d \omega)                                                               \\
	\label{eq:B}
	 & = -\frac{\kappa }{\sigma^2 + 2 \gamma}\  \left[ e^{-\frac{\sigma^2}{2} t}
		                                            \indica{t \geq 0} +  e^{\gamma t} \indica{t < 0} \right] \int_{\R}
	e^{i \omega t } \hat{k}(\gamma + i \omega)g(\d \omega).
\end{align}
%\rc{I would like the denominator to be $\sigma^2/2 +\gamma$ instead of $\sigma^2+\gamma$. This is needed later on. See my next remark.}
Therefore, the Fourier transform of $B$ is given by,
\begin{align*}
	\widehat{B k}(z) & = -\frac{\kappa}{\sigma^2 + 2 \gamma}\int_{\R}
	\frac{\hat{k}(\gamma + i \omega)}{z + \frac{\sigma^2}{2} - i \omega}\  g(\d \omega)
	-\frac{\kappa}{\sigma^2 + 2\gamma}\int_{\R}
	\frac{\hat{k}(\gamma + i \omega)}{\gamma + i \omega -z}\  g(\d \omega),
	\qquad \forall z \in \C
	\ \ \text{with}\ \ \Re(z) = -\lambda.                                                                                \\
	                 & = -\frac{\kappa}{2} \ \int_{\R} \frac{\hat{k}(\gamma + i \omega)}{Q(z- \frac{\beta}{2}, \omega)}\
	g(\d \omega),
\end{align*}
where $\hat{k}(\gamma + i \omega) = \int_0^\infty e^{-(\gamma + i \omega)t }k(t) \d t$.
We now complete the proof of Proposition~\ref{prop:inversibility-I-L}.
\begin{proof}
	Let $\phi \in  L^\infty_\lambda(\R_+)$
	and let $C(\gamma + i \omega_j), j \in I$ be the constants given by Proposition~\ref{prop:inversibility-Laplace-I-L}.
	Set
	\[
		\tilde{\phi}(t) :=  \phi(t)\indica{t \geq 0}
		-\frac{\kappa }{\sigma^2 +2 \gamma}\ \left[ e^{-\frac{\sigma^2}{2} t}
			\indica{t \geq 0} +  e^{\gamma t} \indica{t < 0} \right] \
		\sum_{j \in I} \alpha_j C(\gamma + i \omega_j) e^{i \omega_j t}.
	\]
	We apply Lemma~\ref{lem:invertibility-A} above twice, first with $\theta = \lambda$ and then with $\theta = -\beta -\lambda$. This shows that there exists a Borel measurable function $k$ such that
	\[ \sup_{t \in \R_+} e^{\lambda t} \abs{k(t)} < \infty, \quad \sup_{t \in \R_-} e^{-(\lambda + \beta)t} \abs{k(t)} < \infty \]
	and
	$$
		k(t) = A k(t) + \tilde{\phi}(t), \quad t \in \R.
	$$
	We continue by proving that $k \equiv 0$ on $\R_{-}$.
	Taking the Fourier transform of $k$, we obtain for $\Re(z) = 0$:
	$$
		\hat{k}(z)(1 - \hat{a}(z)) = (1-\hat{a}(z)) \hat{\phi}(z)
		+ \hat{a}(z) \hat{\phi}(z) - \frac{\kappa}{2} \sum_{j \in  I} \frac{\alpha_j C(\gamma + i \omega_j)}{Q(z - \beta/2, \omega_j)},
	$$
	where $Q$ is as in \eqref{eq:Q}.
	Set
	$$
		k(t) = k_1(t) + k_2(t), \qquad k_1(t) = k(t) \indica{t < 0},
		\qquad k_2(t) = k(t)\indica{t \geq 0},
	$$
	Recall that
	\[ a(z) = \frac{\kappa}{2} P(z) = \frac{\kappa}{2} \sum_{j \in I} \frac{\alpha_j}{Q(z-\frac{\beta}{2}, \omega_j)}. \]
	Therefore, we have:
	\begin{equation}
		\label{eq:Phi-on-the-line}
		\hat{k}_1(z) = -\hat{k}_2(z) +  \hat{\phi}(z)
		+ \frac{\kappa/2}{1- \frac{\kappa}{2} P(z)} \sum_{j \in I} \alpha_j \frac{\hat{\phi}(z) - C(\gamma + i \omega_j)}{Q(z - \frac{\beta}{2}, \omega_j)}  ,
		\quad \Re(z) = 0.
	\end{equation}
	By Proposition~\ref{prop:inversibility-Laplace-I-L}, the right-hand-side is holomorphic on $\Re(z) > 0$. Set
	$$
		\Phi(z) = \begin{cases}
			- \hat{k}_2(z) + \hat{\phi}(z) +   \frac{\kappa/2}{1- \frac{\kappa}{2} P(z)} \sum_{j \in I} \alpha_j \frac{\hat{\phi}(z) - C(\gamma + i \omega_j)}{Q(z - \frac{\beta}{2}, \omega_j)} \quad \text{ for } \Re(z) \geq 0, \\
			\hat{k}_1(z), \quad  \text{ for } \Re(z) < 0.
		\end{cases}
	$$
	The function $\Phi$ is continuous on $\Re(z) \geq 0$,
	holomorphic on $\Re(z) > 0$ and bounded.
	Similarly, it is continuous on $\Re(z) \leq 0$,
	holomorphic on $\Re(z) < 0$, and bounded on this
	half space as well. By \eqref{eq:Phi-on-the-line}, $\Phi$ is continuous on $\C$.
	Therefore, by~\cite[th. 16.8]{zbMATH01022658}, $\Phi$ is entire.
	And as it is also bounded, we conclude that $\Phi$ is constant.
	Additionally, as $\Phi(i \omega) \rightarrow 0$ as $\omega \rightarrow \infty$,
	$\Phi \equiv 0$.
	This shows that $k \equiv 0$ on $\R_{-}$.

	Finally, it is shown in Proposition~\ref{prop:inversibility-Laplace-I-L}, $C(\gamma + i \omega_j) = \hat{k}(\gamma + i \omega_j)$.
	Therefore, by the uniqueness of the Laplace transform,
	this shows that $k$ solves
	$$
		k(t) = \kappa L k(t) + \phi(t), \quad t \in \R_+.
	$$
	This completes the proof of the invertibility in $L^\infty_\lambda(\R_+)$ of $I- \kappa L$. The open mapping theorem implies that $(I-\kappa L)^{-1}$ is bounded.
\end{proof}
\begin{remark}
	\label{rem:outline}
	The main idea behind the above proof
	is to write $\kappa L$ as the sum of two operators $A$ and $B$ in such a way that
	the following holds:
	\begin{itemize}
		\item The operator $A$ can be extended ``naturally'' as an operator from $L^1(\R)$ into $L^1(\R)$. The invertibility of $I- A$ is studied using techniques from Fourier analysis.
		\item The operator $B$ from $L^1(\R)$ to $L^1(\R)$ -  when $g$ is the sum of Dirac measures -  is of finite rank. In addition, $I-B$ has an explicit inverse, in the form $I+R$:
		      \[ (I-B)^{-1} \phi = \phi + R \phi. \]
		\item When $k \in L^1(\R_+)$, that is $k \equiv 0$ on $\R_-$, it holds that
		      \[ \kappa Lk = (A+B)k. \]
		      Take $\phi \in L^1(\R_+)$. Our goal is to find $k \in L^1(\R_+)$ such that
		      \[ (I - (A+B)) k  = \phi. \]
		      We rewrite this as $(I - B)(I - (I-B)^{-1}A)k = \phi$, or equivalently (using $I+R = (I-B)^{-1}$):
			      \[ k - (Ak + R (Ak)) = \phi + R \phi, \]
		      that is: $k(t) = (Ak)(t) + \phi(t) + R(\phi + Ak)(t)$.
		\item The operator $R$ above is of rank $2N$
		      and there exist constants $\beta_j \in \C$ such that:
		      $$
			      R (\phi + Ak)(t) =  -\frac{\kappa / 2}{ \frac{\sigma^2}{2} + \gamma}
			      \left[ e^{-\frac{\sigma^2}{2} t} \indica{t \geq 0} +  e^{\gamma t}
				      \indica{t < 0} \right] \sum_{j \in  I} \alpha_j \beta_j e^{i \omega_j t}.
		      $$
		      The ``correct'' constants $\beta_j, j \in I$ are constructed
		      in Proposition~\ref{prop:inversibility-Laplace-I-L}.
	\end{itemize}
\end{remark}}

\section{Estimates}
\label{sec:estimates}
In this section, we gather several technical estimates
and results that are used in the derivation of the
Fr\'echet derivative of the map $\Phi$.  As the optimal
control problem is central in its definition,
we start our analysis with its properties.

\subsection{The optimal control problem}
\label{sec:oc}
To simplify the notation, we denote by $\mathcal{X}_\lambda$ the space of functions of the form:
\[
	h(t,x) = h^1(t) \cos(x) + h^2(t) \sin(x),
\]
where the functions $h^1$ and $h^2$ belong to $L^\infty_\lambda(\R_+)$.
%\rc{I think we should not write $h^1_t$ and $h^2_t$ but we should write $h^1(t)$ and $h^2(t)$. Not only it would be consistent with what we did so far, but it would not tempt the reader to believe that these are derivatives with respect to $t$. I understand that we use time as subscript for stochastic processes like the state $X$ or the Brownian motion, but we did not do that for "functions".}
%\qc{OK I agree it is a good idea. I made the change everywhere, including at the beginning of the text. }
Note that there is a one-to-one correspondence between an element of $L^\infty_\lambda \times L^\infty_\lambda$ and an element of $\cX_\lambda$.
We equip $\cX_\lambda$ with the norm:
\[ \norm{h}_\lambda = \norm{h^1}^\infty_\lambda + \norm{h^2}^\infty_\lambda. \]
For $h \in \mathcal{X}_\lambda$ and $\delta > 0$, we denote by $B_{\mathcal{X}_\lambda}(h;  \delta)$ the open ball $B_{\mathcal{X}_\lambda}(h,  \delta) = \{ k \in \mathcal{X}_\lambda; \quad \norm{k-h}_\lambda < \delta \}$.
Classical parabolic estimates lead to the following existence result of a classical solution of the HJB equation. Its proof is provided in Appendix~\ref{sec:appendix:parabolic} for completeness.
\begin{proposition}
	\label{prop:regular solution HJB}
	Let $\beta, \sigma > 0$. There exists a constant $\eta = \eta(\beta, \sigma) > 0$ such that for all  $\omega \in \R, \lambda > 0$
	and all $f \in B_{\mathcal{X}_\lambda}(0,  \eta)$, the following PDE on $\mathbb{R}_+ \times
		\mathbb{T}$
	\begin{equation}
		\label{eq:HJB}
		v_t - \frac{1}{2} (v_x)^2 + \omega  v_x +  \frac{\sigma^2}{2}v_{xx} + f  = \beta v
	\end{equation}
	has a classical solution $v \in C^{1, 3, 0}(\mathbb{R}_+ \times \mathbb{T} \times \R)$. Moreover, it holds that:
	\begin{enumerate}\itemsep=-1pt
		\item There exists a constant $c_1 = c_1(\beta,\sigma)$ such that for all $t\ge 0$:
		      \[ \sup_{u \geq t} \sup_{x \in \mathbb{T}}  |v_t(u, x, \omega)|\leq c_1( 1 + \abs{\omega})
			      \norm{f}_\lambda e^{-\lambda t}, \quad \omega \in \R. \]
		\item There is a constant $c_2 = c_2(\beta, \sigma)$ such that for all $t \geq 0$ and for all $k \in \{0, 1, 2, 3\}$:
		      \[ \sup_{u \geq t} \sup_{x \in \mathbb{T}, \omega \in \R}  \left| \frac{\partial^k}{\partial x^k } v(u, x, \omega) \right|  \leq c_2 \norm{f}_\lambda e^{-\lambda t}. \]
	\end{enumerate}
\end{proposition}

We continue by defining the control problem
related to \eqref{eq:HJB}.
Let $\eta = \eta(\beta, \sigma)$ be given by
Proposition~\ref{prop:regular solution HJB}. Fix $\lambda > 0$ and
$f \in B_{\mathcal{X}_\lambda}(0,  \eta)$.
Given $\omega \in \R$, $t \geq 0$, $x \in \mathbb{T}$ and $\alpha \in \mathcal{A}_t$,
we define:
$$
	J(t, x, \omega, \alpha) := \E \int_t^\infty{ e^{-\beta(u-t)}
		F(u, X_u, \alpha_u) \d u }, $$
where $F(u, x, a) := f(u, x) + \frac{1}{2} a^2$ and
$X_u = x  + \omega (u-t) + \int_t^u{\alpha_u \d u} + \sigma (B_u-B_t)$.
The following is the classical verification result \cite{FS}.

\begin{proposition}[Verification]
	Let $v(t, x, \omega)$ be a classical solution of \eqref{eq:HJB} as given
	by Proposition~\ref{prop:regular solution HJB}. Then,
	$v(t, x, \omega) \leq J(t, x, \omega, \alpha)$ for every
	$\alpha \in \mathcal{A}_t$, and the
	equality holds if and only if the control is given by
	$$ \alpha_u = -v_x(u, X_u, \omega), \qquad
		\mathbb{P}-\text{a.s., for a.e. } u \geq t.
	$$
\end{proposition}
\begin{proof}
	The proof follows a standard \emph{verification} argument. First, we have  $-\frac{1}{2} v^2_x = \inf_{\alpha \in
			\mathbb{R}} \alpha v_x + \frac{1}{2} \alpha^2$.
	So, if we define $\mathcal{L}_\alpha v := v_t + \alpha v_x + \omega v_x + \frac{\sigma^2}{2} v_{xx}$, we have:
	\begin{equation}
		\label{eq:property HJB}
		\forall \alpha \in \mathbb{R}, \quad  \beta v \leq \mathcal{L}_\alpha v + F(t, x, \alpha).
	\end{equation}
	By Ito's lemma, we have for all $\theta \geq t$:
	\begin{align*}
		e^{-\beta(\theta-t)} v(\theta, X_\theta, \omega) = v(t, x, \omega)  + & \int_t^\theta{ e^{-\beta(u-t)}[-\beta v(u, X_u, \omega) +
					                                                                               \mathcal{L}_{\alpha_u}v (u, X_u, \omega) ] \d u}                 \\
		                                                                      & \quad + \int_t^\theta{ e^{-\beta(u-t)} v_x(u, X_u, \omega) \sigma \d B_u}.
	\end{align*}
	Because $v_x$ is bounded, the last term is a true martingale. So \eqref{eq:property HJB} yields:
	\[ \E e^{-\beta(\theta-t)} v(\theta, X_\theta, \omega) \geq v(t, x, \omega) - \E \int_t^\theta{ e^{-\beta(u-t)}F(u, X_u, \alpha) \d u}.
	\]
	We now let $\theta \rightarrow \infty$ to obtain the stated inequality.
\end{proof}
\begin{remark}
	Note that this result implies the uniqueness of the solution of \eqref{eq:HJB} because if $v$ solves
	\eqref{eq:HJB} and satisfies the growth condition:
	\[ \lim_{t \rightarrow \infty} \sup_{x \in \mathbb{T}} |v(t, x, \omega)| e^{-\beta t} = 0, \] we necessarily have
	$v(t, x, \omega) = \inf_{\alpha \in \mathcal{A}_t} J(t, x, \omega, \alpha)$.
	In addition, $v_x(t, \cdot, \omega)$ is Lipschitz on $\mathbb{T}$, uniformly with respect to $t$, so the SDE
	\[ X^{x, t, *}_u = x + \omega(u-t) - \int_t^u{ v_x(\theta, X^{x, t, *}_\theta, \omega) \d \theta} + \sigma(B_u - B_t)\]
	has a unique path-wise solution, defined for all $u \geq t$.
	This shows that the minimum of $J$ is attained for $\alpha_u = -v_x(u, X^{x, t, *}_u, \omega)$.
\end{remark}

\subsection{Estimates using Pontryagin's maximum principle}
\label{ss.prntryagin}
This technical subsection provides the central
estimates for the sensitivity analysis of the next subsection.

Let $\eta = \eta(\beta, \sigma)$ be given by Proposition~\ref{prop:regular solution HJB}. Fix $\lambda > 0$
and $f,h \in B_{\mathcal{X}_\lambda}(0,  \eta/2)$.
Given $\theta > 0$, we denote by $h^{[\theta]}$ the following function:
$$
	h^{[\theta]}(u, x) := h(u, x)\mathbbm{1}_{u \leq \theta},
	\qquad \text{for}\  u \geq 0,  x \in \mathbb{T}.
$$
We let $v^{[\theta]}(t, x, \omega)$ be the value function associated with the cost
$f+h^{[\theta]} + \frac{\alpha^2}{2}$.
In other words, $v^{[\theta]}$ solves on $\mathbb{R}_+ \times \mathbb{T} \times \R$:
\[ v^{[\theta]}_t - \frac{1}{2} (v^{[\theta]}_x)^2 + \omega v^{[\theta]}_x +  \frac{\sigma^2}{2} v^{[\theta]}_{xx} + f + h^{[\theta]} =
	\beta
	v^{[\theta]}. \]
Note that this construction implies that:
\[ v^{[t]}(t, x, \omega) = v^{[0]}(t, x, \omega), \quad \forall t \geq 0. \]
To simplify the notation, we write $v^f(t, x, \omega) := v^{[0]}(t, x, \omega)$ and $v^{f+h}(t, x, \omega) :=  v^{[\infty]}(t, x, \omega)$.
We denote by
$(X^{t, x, [\theta]}_u)$ the
solution of the SDE
\begin{equation} \forall u \geq t, \quad X^{x, t, [\theta]}_u = x + \omega(u-t) - \int_t^u{ v^{[\theta]}_x(s, X^{x, t,
						[\theta]}_s, \omega) \d s} + \sigma (B_u - B_t),
	\label{eq:def of X xt[theta]}
\end{equation}
and write $Y^{x, t, [\theta]}_u := v^{[\theta]}_x(u, X^{x, t, [\theta]}_u, \omega)$ and $Z^{x, t, [\theta]}_u = \sigma
	v^{[\theta]}_{xx}(u, X^{x, t, [\theta]}_u, \omega)$.
Ito's formula yields
\begin{align*}
	\d Y^{x, t, [\theta]}_u   & = \beta Y^{x, t, [\theta]}_u \d u  -(f+h^{[\theta]})_x(u, X^{x, t, [\theta]}_u) \d u + Z^{x, t, [\theta]}_u \d B_u \\
	Y^{x, t, [\theta]}_\theta & = v^{[\theta]}_x(\theta, X^{x, t, [\theta]}_\theta, \omega).
\end{align*}
This BSDE corresponds to the \textit{adjoint equation} of the Pontryagin maximum principle.
We have the following estimates.
\begin{proposition}
	\label{prop:estimates theta versus theta+epsilon}
	Let $\lambda, \theta > 0$ be fixed and $h \in B_{\mathcal{X}_\lambda}(0,  \eta/2)$. There exist two constants
	$\delta = \delta(\beta, \sigma, \lambda) > 0$ and $C = C(\beta, \lambda)$ such that for all $f \in B_{\mathcal{X}_\lambda}(0,  \delta)$, for all $t \in [0, \theta]$, $\omega \in \R$:
	$$ | X^{t, x, [\theta]}_u - X^{t, x, [0]}_u|
		\leq C \int_t^\theta{ \norm{h(v)}_\infty \d v},
		\qquad \forall u \in [t, \theta],
	$$
	where for  $v \geq 0$, $\norm{h(v)}_\infty = \sup_{x \in \T} \abs{h(v, x)}$.
\end{proposition}
\begin{proof}
	We write $X_u := X^{x, t, [\theta]}_u - X^{x, t, [0]}_u$, $Y_u := Y^{x, t, [\theta]}_u - Y^{x,
				t, [0]}_u$, $Z_u := Z^{x, t, [\theta]}_u - Z^{x, t, [0]}_u$, and:
	\begin{align*}  \xi^1_u      & := h_x(u, X^{x, t, [\theta]}_u),                                                                \\
	                \xi^2_u      & := \frac{1}{X_u} \left[ f_x(u, X^{x, t, [\theta]}_u) -  f_x(u,
		                                                X^{x, t, [0]}_u) \right],                                  \\
	                \xi^3_\theta & :=  \frac{1}{X_\theta} \left[ v^f_x(\theta, X^{x, t, [\theta]}_\theta, \omega) -  v^f_x(\theta,
		                                                      X^{x, t, [0]}_\theta, \omega) \right].
	\end{align*}
	Take $\delta \leq \beta^2/8$. We have:
	\[ |\xi^1_u| \leq  \norm{h(u)}_\infty \quad \text{ and } \quad |\xi^2_u| \leq \frac{\beta^2}{8} e^{-\lambda u}. \]
	In addition,  Proposition~\ref{prop:regular solution HJB} ensures that provided that $\norm{f}_\lambda$ is small enough:
	\[ |\xi^3_{\theta}| \leq c_2 \norm{f}_\lambda e^{-\lambda \theta} \leq \frac{\beta}{8} e^{-\lambda \theta}. \]
	%{\color{red} Here I would state the bounds on $|\xi^1_u|$ and $|\xi^2_u|$.}
	Using the notations above, it holds that for all $u \in [t, \theta]$:
	\begin{align*}
		\d X_u & = - Y_u \d u, \quad X_t = 0,                                                                                 \\
		\d Y_u & = \beta Y_u \d u -  \xi^1_u \d u -  \xi^2_u X_u \d u + Z_u \d B_u, \quad Y_{\theta} = \xi^3_\theta X_\theta.
	\end{align*}
	We view this system as a linear FBSDE driven by the random coefficients $\xi^1$, $\xi^2$ and $\xi^3_\theta$.
	Following \cite[Ch. 2.6]{Car}, we look for solutions of the form $Y_u = P_u X_u+ p_u$ for
	adapted processes $P$ and $p$. We find that $P$ and $p$ solve the following equations:
	\begin{align}
		\label{eq:Pt BSDE}
		\d P_u & = (-\xi^2_u + \beta P_u + (P_u)^2) \d u + \Lambda_u \d B_u, \quad P_\theta = \xi^3_\theta \\
		\d p_u & = (\beta p_u + P_u p_u - \xi^1_u) \d u + \tilde{\Lambda} _u \d B_u, \quad p_\theta = 0,
		\label{eq:pt BSDE}
	\end{align}
	for some adapted processes $\Lambda, \tilde{\Lambda}$ satisfying $Z_u = \Lambda_u X_u + \tilde{\Lambda}_u$.
	%\rc{I would replace the process $\lambda$ by $\tilde\Lambda$ to avoid confusion with the constant $\lambda$. This should not be a big change because if I am not mistaken, it does not appear below.}
	The first equation only involves $(P_t)$. It is a quadratic BSDE. It is well posed for the following reasons:
	\begin{enumerate}
		\item The coefficient of $P^2$ being $1$ (positive);
		\item $\xi^2_u$ is uniformly bounded in $u$ and $\omega$;
		\item the terminal condition $\xi^3_\theta$ is bounded.
	\end{enumerate}
	Under these conditions, standard arguments can be used to prove a priori upper and lower bounds for the solution, and once we know the solution remains in a bounded interval, the quadratic term can be controlled as it becomes impotent. From this point, any existence and uniqueness result can be used.

	\textbf{Claim 1}.
	Since $|\xi^2_u| \leq \frac{\beta^2}{8} e^{-\lambda
				u}$ and $|\xi^3_\theta| \leq \frac{\beta}{8} e^{-\lambda
				\theta}$, the solution of the BSDE~\eqref{eq:Pt BSDE} satisfies:
	\[ \mathbb{P}-a.s., \quad |P_u| \leq \frac{\beta}{2} e^{-\lambda u}. \]
	Indeed, define:
	\[
		P^*_u := \text{esssup}_{\Omega} \sup_{s \in [u, \theta]} |P_s|
	\]
	and $\tau = \sup\{ u \in [t, \theta] ~|~  \frac{(P^*_u)^2}{\beta} > \frac{|P^*_u|}{2} \}$.
	From \eqref{eq:Pt BSDE}, we have
	\[ P_u = \E \left[e^{-\beta (\theta-u)} \xi^3_\theta ~|~ \mathcal{F}_u \right] +  \E \left[ \int_u^\theta{ e^{-\beta
						(s-u)}(\xi^2_s -
				(P_s)^2) \d s} ~|~ \mathcal{F}_u
			\right],  \]
	and so we deduce that on $[\tau, \theta]$
	\[ |P^*_u| \leq  \frac{\beta}{8} e^{-\lambda \theta} +   \frac{\beta}{8} e^{-\lambda u} + \frac{1}{2} |P^*_u|, \]
	and so
	\[ |P^*_u| \leq \frac{\beta}{2} e^{-\lambda u}. \]
	We deduce from this bound that $\tau = t$ and so the claim is proven. %{\color{red} Not sure I got that}

	\vskip 1pt\noindent
	\textbf{Claim 2}. It holds that
	\[ \forall u \in [t, \theta], \quad |p_u| \leq  \int_u^\theta{ e^{-\frac{\beta}{2}(s-u)} \norm{h(s)}_\infty \d s}. \]
	Indeed, from \eqref{eq:pt BSDE}, we have
	\[ p_u =  \E \left[ \int_u^\theta { \Gamma^s_u \xi^1_s \d s} ~|~ \mathcal{F}_u \right],\]
	where $\Gamma^s_u := \exp \left(-\int_u^s{ P_v \d v} \right) e^{-\beta (s-u)}$.
	%{\color{red}I am not sure if the formula for $\Gamma^s_u$ is correct. I am not sure of the sign of $P_v$, I would 
	%put a $-$. In any case, I am sure that the upper bound is correct in the sense that $|\Gamma^s_u| \leq \exp \left( 
	%\int_u^s{ |P_v| dv} \right) e^{-\beta (s-u)}$.}
	By Claim 1, we have $|\Gamma^s_u| \leq e^{-\frac{\beta}{2}(s-u)}$.
	We deduce that
	\[ |p_u| \leq \int_u^\theta{ e^{-\frac{\beta}{2}(s-u)} \text{esssup}_{\Omega} |\xi^1_s| \d s}. \]
	This concludes the proof of the claim.\\
	\textbf{Claim 3}. There exists a constant $C > 0$ such that for all $u \in [t, \theta]$:
	\[ |X_u| \leq C \int_t^\theta{ \norm{h(v)}_\infty \d v}. \]
	%\rc{$\int_t^\theta{ \norm{h(v)}_\infty \d v}$ for me.}
	Recall that $Y_u= P_u X_u + p_u$ and $\d X_u= - Y_u \d u$. Using that $X_t = 0$, we deduce that
	\[ X_u = \int_t^u{ e^{-(A_u - A_s)} p_s \d s}, \quad \text{ where } \quad  A_u-A_s :=  \int_s^u{ P_v \d v}.  \]
	By Claim 1, we have $|A_u - A_s| \leq \frac{\beta}{2 \lambda}$ and so by Claim 2 we have
	\begin{align*}
		|X_u| & \leq e^{\frac{\beta}{2 \lambda}} \int_t^u{\int_s^\theta {e^{-\frac{\beta}{2}(v-s)} \norm{h(v)} \d v} \d s}
		\\ &	\leq \frac{2e^{\frac{\beta}{2 \lambda}} }{\beta}   \int_t^\theta{ \norm{h(v)} \d v}.
	\end{align*}
	%\rc{Again, I would replace $\norm{h_v}$ by $\norm{h(v)}_\infty$ for me.}
	This concludes the proof.
\end{proof}

\subsection{Sensitivity Analysis}
\label{sec:sensitivity analysis}

Let $\lambda > 0$.
Let $\eta = \eta(\beta, \sigma)$ be given by Proposition~\ref{prop:regular solution HJB} and $\delta = \delta(\beta,
	\sigma, \lambda)$ be given by Proposition~\ref{prop:estimates theta versus theta+epsilon}.
We fix $f,h \in B_{\mathcal{X}_\lambda}(0,  \frac{\eta \wedge \delta}{2})$. This assumption will be in force
throughout the end of this section.
We first study the value function and then the SDE
governing the optimal state process.

Recall that for $\theta > 0$, $h^{[\theta]}(t, x) = h(t, x)\mathbbm{1}_{t \leq
		\theta}$ and $v^{[\theta]}(t, x, \omega)$ is the value function associated to the cost $f(t, x) + h^{[\theta]}(t, x) +
	\alpha^2/2$. We write:
\[ v^{[\theta]}(t, x, \omega) = \inf_{\alpha \in \cA_t} J^{[\theta]}(t, x, \omega, \alpha). \]
Recall also that we use the notations $v^{f+h}(t, x, \omega) = v^{[\infty]}(t, x, \omega)$ and $v^f(t, x, \omega) = v^{[t]}(t, x, \omega)$.
The main result of this section is the following:
\begin{proposition}
	\label{prop:sensitivity estimates value function}
	Let $X^{x, t, [\theta]}_u$ be the solution of SDE \eqref{eq:def of X xt[theta]}. It holds that
	\[ v^{f+h}(t, x, \omega) - v^f(t, x, \omega) = \int_t^\infty{ e^{-\beta(\theta - t)} \E h(\theta, X^{x, t, [\theta]}_\theta)
			\d \theta}.
	\]
\end{proposition}
\begin{proof}
	Let $\epsilon  \in (0, 1)$, $x \in \mathbb{T}$ and $\theta > t$ be fixed.
	Let $\alpha \in \mathcal{A}_t$ be defined by
	\[ \alpha_u := -v^{[\theta]}_x(u, X^{x, t, [\theta]}_u, \omega). \] We have
	\begin{align*}
		v^{[\theta+\epsilon]}(t, x, \omega) - v^{[\theta]}(t, x, \omega) & \leq J^{[\theta+\epsilon]}(t, x, \omega, \alpha) -
		v^{[\theta]}(t, x, \omega)                                                                                                                                      \\
		                                                                 & = \E \int_t^\infty{ e^{-\beta (u-t)} \left[ h^{[\theta+\epsilon]}(u, X^{x, t, [\theta]}_u) -
				                                                                               h^{[\theta]}(u, X^{x, t, [\theta]}_u) \right] \d u} \\
		                                                                 & = \int_\theta^{\theta+\epsilon}{ e^{-\beta (u-t)}  \E h(u, X^{x, t, [\theta]}_u) \d u  }.
	\end{align*}
	The function $u \mapsto \E h(u, X^{x, t, [\theta]}_u)$ is continuous at $u = \theta$, therefore
	\[ \limsup_{\epsilon \downarrow 0}  \frac{v^{[\theta+\epsilon]}(t, x, \omega) - v^{[\theta]}(t, x, \omega)}{\epsilon} \leq
		e^{-\beta(\theta-t)} \E h(\theta, X^{x, t, [\theta]}_\theta). \]
	Similarly, we have, choosing this time $\alpha_u := -v^{[\theta+\epsilon]}_x(u, X^{x, t, [\theta+\epsilon]}_u, \omega)$:
	\begin{align*}
		v^{[\theta+\epsilon]}(t, x, \omega) - v^{[\theta]}(t, x, \omega) & \geq v^{[\theta+\epsilon]}(t, x, \omega) - J^{[\theta]}(\alpha, t, x, \omega)                   \\
		                                                                 & = \int_\theta^{\theta+\epsilon} {e^{-\beta(u-t)} \E h(u, X^{x, t, [\theta+\epsilon]}_u) \d u }.
	\end{align*}
	Using Proposition~\ref{prop:estimates theta versus theta+epsilon} (with $\theta$ replaced by $\theta+\epsilon$,
	$h$ replaced by $h =
		h^{[\theta+\epsilon]} -
		h^{[\theta]}$ and $f$ replaced by $f = f + h^{[\theta]}$), we find that there exists a constant $C(\lambda, \beta,
		h)$ such
	that
	\[ |X^{x, t, [\theta+\epsilon]}_u - X^{x, t, [\theta]}_u| \leq C(\lambda, \beta, h) \epsilon. \]
	The function $h$ being Lipschitz with respect to $x$, we deduce that
	\[ \liminf_{\epsilon \downarrow 0}  \frac{v^{[\theta+\epsilon]}(t, x, \omega) - v^{[\theta]}(t, x, \omega)}{\epsilon} \ge
		e^{-\beta(\theta-t)} \E h(\theta, X^{x, t, [\theta]}_\theta). \]
	This proves that $\theta \mapsto v^{[\theta]}(t, x)$ is right-differentiable. Left-differentiability is proved
	similarly. We deduce that $\theta \mapsto v^{[\theta]}(t, x, \omega)$ is differentiable everywhere on $(t, \infty)$ with
	\[ \frac{\d}{\d \theta} v^{[\theta]}(t, x, \omega) = e^{-\beta(\theta-t)} \E h(\theta, X^{x, t, [\theta]}_\theta). \]
	Integrating this equality  from $\theta = t$ to $\theta = \infty$ concludes the proof.
	%$\rc{I thought $v^{[\theta]}=v^f$ for $\theta =0$, not $\theta=t$. Am I missing something?} \qc{OK, you are right. I changed a little bit the presentation - but in any case we have $v^{[0]}(t, x,  \omega) = v^{[t]}(t, x, \omega), \forall t$, as in both optimization problem the additional cost $h$ is always equal to zero.}
\end{proof}

The following results list the main estimates that will be used later on.
\begin{corollary}
	\label{cor:estimate vxh minus vx}
	It holds that
	\begin{equation}
		\label{eq:v_xf+h - vxf}
		v^{f+h}_x(t, x, \omega) -  v^f_x(t, x, \omega) =  \int_t^\infty{ e^{-\beta(\theta - t)} \E \left[ h_x(\theta, X^{x, t, [\theta]}_\theta)
				\frac{\d}{\d x} X^{x, t, [\theta]}_\theta \right] \d \theta},
	\end{equation}
	with
	\begin{equation}
		\label{eq:dxX}
		\frac{\d}{\d x} X^{x, t, [\theta]}_\theta = \exp\left( - \int_t^\theta v^{[\theta]}_{xx}(s, X^{x, t, [\theta]}_s, \omega) \d s
		\right).
	\end{equation}
	Moreover, by Proposition~\ref{prop:regular solution HJB}, for all $f,h \in
		B_{\mathcal{X}_\lambda}(0, \eta/2)$,
	\[  \forall \theta \geq 0, \quad   \sup_{t \geq 0} \sup_{x \in \mathbb{T}, \omega \in \R }|v^{[\theta]}_{xx}(t, x, \omega)| e^{\lambda t}
		\leq C(\eta). \]
	So, there exists a constant $C(\eta, \lambda)$ such that for all $f,h \in
		B_{\mathcal{X}_\lambda}(0, \eta/2)$:
	\[  \sup_{\theta \geq 0} \sup_{t \geq 0}\sup_{x \in \mathbb{T}, \omega \in \R} \left| \frac{d}{dx} X^{x, t, [\theta]}_\theta \right|
		\leq C(\eta,\lambda)\]
	and consequently,
	\begin{equation}
		\label{eq: estimate vx f+h - vx f}
		| v^{f+h}_x(t, x) -  v^f_x(t, x)| \leq \frac{C(\eta,\lambda)}{\beta} \norm{h}_\lambda
		e^{-\lambda t}.
	\end{equation}
\end{corollary}
\vspace{5pt}

We continue with the analysis of the  forward SDE.
Denote by $v^{f}(t, x, \omega)$ the value function associated to the cost $f(t, x)+\frac{1}{2} \alpha^2$. Given $\mu \in
	\mathcal{P}(\mathbb{T})$, let
$(X^{\mu,\omega,f}_{t, s})_{t \geq s}$ be the solution of
\[ \d X^{\mu,\omega,f}_{t,s} =  \omega \d t - v^f_x(t, X^{\mu,\omega,f}_{t, s}, \omega) \d t+\sigma \d B_t, \quad \text { with } \quad X^{\mu,\omega,f}_{s,s} \sim
	\mu.
\]
We now follow \cite[Prop. 2.11]{Cor} to obtain an estimate of the state process.
\begin{proposition}
	\label{prop:sensitivity formula SDE}
	For any $g \in C^2(\mathbb{T})$ and $\mu \in \cP(\T)$, for all $t \geq s$:
	\begin{align*}
		\E [g(X^{\mu,\omega,f+h}_{t, s})- g(X^{\mu,\omega,f}_{t, s})]
		= - \int_s^t{ \int_{\mathbb{T}}
				           \left[ \frac{\d}{\d y} \E g(X^{\delta_y, \omega, f}_{t, \theta})  \right]
				           (v^{f+h}_x  - v^f_x)(\theta, y, \omega) \mathcal{L}(X^{\mu,\omega,f+h}_{\theta, s})
				           (\d y) \d \theta}.
	\end{align*}
\end{proposition}

\section{Fr\'echet derivative}
\label{sec:Frechet}

Let $\lambda > 0$.
Recall the function $\Phi$ of subsection~\ref{ss.reformulation}. The above notations imply that
\[ \E \cos(X^{\omega, *}_t) = \E \cos(X^{\nu^\omega, \omega, h}_{t, 0}). \]
That is, with a slight abuse of notations, we can rewrite $\Phi$ as:
\[
	\Phi(\nu, h) = - h + \kappa \begin{pmatrix}  \int \E \cos(X^{\nu^\omega, \omega,h}_{t, 0}) g(\d \omega) \\ \int{\E \sin(X^{\nu^\omega, \omega,h}_t) g(\d \omega)}  \end{pmatrix}_{t \geq 0},
	\qquad  h \in \cX_\lambda\;\; \nu \in \mathcal{P}(\mathbb{T}\times \mathbb{R}),
\]
where the second marginal of $\nu$ is equal to $g$, namely $\nu(\d x,\d \omega)=\nu^\omega(\d x)g(\d \omega)$. Recall the pseudo distance $\dist$ defined by \eqref{eq:def-dist}.
Let $\epsilon > 0$. We assume that
\[ g(\d \omega) a.e. , \quad \dist(\nu^\omega) \leq \epsilon. \]
\begin{lemma}
	\label{lem:into}
	For $\lambda \in (0, \sigma^2/2)$,
	there is $\eta$ small enough such that
	$\Phi(\nu, h) \in L^\infty_\lambda \times L^\infty_\lambda$
	whenever  $\norm{h}^\infty_\lambda < \eta$.
\end{lemma}
\begin{proof}
	Proposition~\ref{prop:sensitivity formula SDE} with $g = \cos$ and $f = 0$ gives:
	\[ 	\E \cos(X^{\nu^\omega,\omega, h}_{t, 0}) = \E \cos(X^{\nu^\omega, \omega, 0}_{t, 0}) -\int_0^t{ \int_{\mathbb{T}}  \left[  \frac{\d}{\d y} \E
					\cos(X^{\delta_y, \omega, 0}_{t, \theta})
					\right](v^h_x
				- v^0_x)(\theta, y, \omega)
				\mathcal{L}(X^{\nu^\omega,\omega, h}_{\theta, 0}) (\d y) \d \theta}.   \]
	The value function associated with the cost $\frac{1}{2} \alpha^2$ is equal to zero. We deduce that $v^0_x =
		0$ and $X^{\delta_y, \omega, 0}_{t, \theta} = y + \omega(t-\theta) + \sigma (B_t - B_\theta)$. So
	\[ \E
		\cos(X^{\delta_y, \omega, 0}_{t, \theta})  = \cos(y + \omega (t-\theta)) e^{-\frac{\sigma^2}{2} (t-\theta)}, \]
	and
	\begin{align}
		\label{eq:expectation cos}
		\E \cos(X^{\nu^\omega, \omega, h}_{t, 0}) & = \int_{\T} \cos(y + \omega t) \nu^\omega(\d y) e^{-\frac{\sigma^2}{2} t}                        \\
		                                          & \quad +  \int_0^t{ e^{-\frac{\sigma^2}{2}
						                                                            (t-\theta)} \int_{\mathbb{T}} \sin(y + \omega(t-\theta))
				                                                            v^h_x(\theta, y, \omega)
				                                                            \mathcal{L}(X^{\nu^\omega,\omega, h}_{\theta, 0}) (\d y) \d \theta}. \nonumber
	\end{align}
	Similarly, we have
	\begin{align}
		\label{eq:expectation sin}
		\E \sin(X^{\nu,\omega, h}_{t,0}) & = \int_{\T} \sin(y + \omega t) \nu(\d y) e^{-\frac{\sigma^2}{2} t}                        \\
		                                 & \quad -  \int_0^t{ e^{-\frac{\sigma^2}{2}
						                                                   (t-\theta)} \int_{\mathbb{T}} \cos(y + \omega(t-\theta))
				                                                   v^h_x(\theta, y, \omega)
				                                                   \mathcal{L}(X^{\nu,\omega, h}_{\theta, 0}) (\d y) \d \theta}. \nonumber
	\end{align}
	By Corollary~\ref{cor:estimate vxh minus vx}, we have:
	\[ \forall y \in \mathbb{T}, \forall t \geq 0, \quad  |v_x^h(t, y)| \leq  C(\beta, \lambda, \delta)
		\norm{h}_\lambda e^{-\lambda t}. \]
	Therefore:
	\begin{align} \left| \E \cos(X^{\nu, \omega, h}_{t,0}) \right| & \leq e^{-\frac{\sigma^2}{2} t} \dist(\nu^\omega) + C(\beta, \lambda, \delta) \norm{h}_\lambda
	              \int_0^t{
			                     e^{-\frac{\sigma^2}{2}
					                     (t-\theta)}  e^{-\lambda \theta} \d
			                     \theta}  \nonumber                                                                                                                            \\
	                                                               & \leq e^{-\frac{\sigma^2}{2} t} \dist(\nu^\omega)  + \frac{C(\beta, \lambda, \delta)}{ \frac{\sigma^2}{2} - \lambda } \norm{h}_\lambda
	              e^{-\lambda t}. \label{eq:estimate 06/05}
	\end{align}
	A similar estimate holds for $\left| \E \sin(X^{\nu, \omega, h}_{t, 0}) \right| $. This concludes the proof.
\end{proof}

The following gradient estimate is used in the subsequent analysis.

\begin{lemma}[Gradient estimate]
	\label{lem:grad estimate}
	Fix $\lambda \in (0, \sigma^2/2)$, and  either $g=\cos$ or  $g=\sin$.
	There is a constant $C(\eta, \lambda)$ such that for all $h \in B_{\mathcal{X}_\lambda}(0,  \eta)$ and $x \in \mathbb{T}$,
	$$
		\frac{\d }{\d x} \E g(X^{\delta_x, \omega, h}_{t, s}) \leq C(\eta, \lambda)  \left( e^{-\frac{\sigma^2}{2} (t-s)} +e^{-\lambda t} \right),
		\qquad \forall t \geq s \geq 0.
	$$
	A similar estimate holds for $\frac{\d^2 }{\d x^2} \E g(X^{\delta_x, \omega, h}_{t, s})$.
\end{lemma}
\begin{proof}
	We use Proposition~\ref{prop:sensitivity formula SDE} with $g = \cos$ and
	$\nu = \delta_x$ to arrive at
	$$
		\E \cos(X^{\delta_x, \omega, h}_{t, s})  = \cos(x + \omega(t-s))
		e^{-\frac{\sigma^2}{2} (t-s)} +  \int_s^t e^{-\frac{\sigma^2}{2}
				(t-\theta)} \E \sin(X^{\delta_x, \omega, h}_{\theta, s}
		+ \omega(t-\theta)) v^h_x(\theta, X^{\delta_x,\omega,h}_{\theta, s}) \d \theta.
	$$
	Therefore,
	\begin{align} \frac{\d}{\d x} \E \cos(X^{\delta_x, \omega, h}_{t, s}) & =  -\sin(x + \omega(t-s)) e^{-\frac{\sigma^2}{2} (t-s)} \label{eq:repr-cos_x}                                                                            \\
	                                                                      & \quad +  \int_s^t{
			                                                                                        e^{-\frac{\sigma^2}{2}
					                                                                                        (t-\theta)} \E \cos(X^{\delta_x,\omega,h}_{\theta, s}) \left( \frac{\d}{\d x}X^{\delta_x,\omega,h}_{\theta, s} \right) v^h_x(\theta,
			                                                                                        X^{\delta_x,\omega,h}_{\theta,s}) \d \theta} \nonumber \\
	                                                                      & \quad +\int_s^t{ 	e^{-\frac{\sigma^2}{2}
					                                                                                      (t-\theta)}    \E \sin(X^{\delta_x,s,h}_\theta)v^h_{xx}(\theta,
			                                                                                      X^{\delta_x,\omega,h}_{\theta, s})  \frac{\d}{\d x}X^{\delta_x,\omega,h}_{\theta, s}  \d \theta}.  \nonumber
	\end{align}
	In view of Proposition~\ref{prop:regular solution HJB},
	$|v^h_{x}(t, x)|  + |v^h_{xx}(t, x)| \leq C(\eta, \sigma, \lambda) e^{-\lambda t}$,
	and by Corollary~\ref{cor:estimate vxh minus vx} we have
	$\left| \frac{\d}{\d x}X^{\delta_x,\omega,h}_{\theta, s} \right|
		\leq  C(\eta, \sigma, \lambda)$. These  yield the first estimate.
	By differentiating \eqref{eq:repr-cos_x} with respect to $x$, we deduce the similar estimate for  $\frac{\d^2 }{\d x^2} \E g(X^{\delta_x, \omega, h}_{t, s})$.
\end{proof}
\begin{lemma}
	\label{lem:diff-vf}
	For $\lambda \in (0, \sigma^2/2)$,
	the map $B_{\mathcal{X}_\lambda}(0, \eta) \ni
		f \mapsto  v^f_x(t, x, \omega) \in \mathbb{R}$ is
	Fr\'echet differentiable, and
	\[  D_f v^f_x(t, x, \omega) \cdot h =  \int_t^\infty{ e^{-\beta(\theta - t)} \E \left[ h_x(\theta, X^{\delta_x, \omega, f}_{\theta, t})
				\exp\left( -\int_t^\theta{ v^f_{xx}(s, X^{\delta_x, \omega, f}_{s, t}, \omega) \d s} \right) \right]
			\d \theta}.  \]
\end{lemma}
\begin{proof}
	We prove that there exists a constant $C = C(\beta, \sigma, \eta, \lambda)$ such that for all $f, h \in
		B_{\mathcal{X}_\lambda}(0, \eta/2)$:
	\begin{equation}
		\label{eq:estimate vx f+h - vx f - Dvxf}
		\left| v^{f+h}_x(t, x, \omega) - v^{f}_x(t, x, \omega) - D_f v^f_x(t, x, \omega) \cdot h \right| \leq
		C  \norm{h}^2_\lambda e^{-2 \lambda t}.
	\end{equation}
	The proof directly follows from~\eqref{eq:v_xf+h - vxf} and the previous estimates.
	Indeed, by Prop.~\ref{prop:estimates theta versus theta+epsilon},
	$$
		|h_x(\theta, X^{x, t, [\theta]}_\theta) -h_x(\theta, X^{x, t, [0]}_\theta) |
		\leq  \norm{h}_\lambda e^{-\lambda \theta}
		\left| X^{x,t, [\theta]}_{\theta} - X^{x,t, [0]}_{\theta} \right|
		\leq   C \norm{h}^2_\lambda e^{-2 \lambda t}    (\theta-t).
	$$
	Similarly, \eqref{eq:v_xf+h - vxf} yields
	$$
		|v^{[\theta]}_{xx}(t, x) - v^{[0]}_{xx}(t, x)| \leq C
		\norm{h}_\lambda e^{-\lambda t},
	$$
	and by  \eqref{eq:dxX},
	$$
		\left| \frac{d}{dx} X^{x, t, [\theta]}_\theta - \frac{d}{dx} X^{x, t, [0]}_{\theta}
		\right| \leq C 	\norm{h}_\lambda e^{-\lambda t}(\theta -t).
	$$
	We complete the proof by using these estimates
	together with \eqref{eq:v_xf+h - vxf}.
\end{proof}

\begin{corollary}
	\label{cor:derivative-v_x_near_0}
	For $h(t,x) = h^1(t) \cos(x) + h^2(t) \sin(x) \in \mathcal{X}_\lambda$,
	$$
		D_f v^{0}_x(t, x, \omega) \cdot h
		=  (- K^2_\omega * h^1 + K^1_\omega * h^2)(t) \cos(x) -
		(K^1_\omega * h^1 + K^2_\omega * h^2)(t) \sin(x),
	$$
	where $*$ denotes the convolution operator and
	$$
		K^1_\omega(t) := e^{(\beta + \frac{\sigma^2}{2}) t} \cos(\omega t)
		\mathbbm{1}_{(-\infty,0]}(t), \qquad
		K^2_\omega(t) := e^{(\beta + \frac{\sigma^2}{2}) t} \sin(\omega t)
		\mathbbm{1}_{(-\infty,0]}(t),
		\qquad  t \in \mathbb{R}.
	$$
\end{corollary}
\begin{proof}
	We have $X^{\delta_x, \omega, 0}_{\theta, t} = x + \omega(\theta-t) + \sigma B_{\theta-t}$. Therefore,
	\[ \E \cos(X^{\delta_x, \omega, 0}_{\theta, t}) = \left[ \cos(x) \cos(\omega(\theta-t)) - \sin(x) \sin(\omega(\theta-t)) \right] e^{-\frac{\sigma^2}{2} (\theta-t)}. \]
	A similar formula holds for $\E \sin(X^{\delta_x, \omega, 0}_{\theta, t})$.
	Combining these we arrive at
	\begin{align*} \E h_x(\theta, X^{\delta_x, \omega, 0}_{\theta, t}) & = \left[ -h^1(\theta) \sin(x) \cos(\omega(\theta-t)) -h^1(\theta) \cos(x) \sin(\omega (\theta-t)) \right. \\
		                                                                       & \quad +  \left. h^2(\theta) \cos(x) \cos(\omega (\theta-t)) - h^2(\theta) \sin(x) \sin(\omega(\theta-t)) \right] e^{-\frac{\sigma^2}{2} (\theta-t)}.
	\end{align*}
	We deduce the result by applying Lemma~\ref{lem:diff-vf}.
\end{proof}

We are now ready to prove the central result of this section.

\begin{proposition}
	\label{prop:Phi is Fr\'echet differentiable}
	The function
	$\Phi(\nu, \cdot): \mathcal{B}_{\mathcal{X}_\lambda}(0,\rho) \rightarrow \mathcal{X}_\lambda$ is Fr\'echet differentiable and
	$$
		(D_f \Phi(\nu, f) \cdot h)(t, x)  =
		-h(t, x) + L^{1, \nu}_t(h) \cos(x) + L^{2, \nu}_t(h) \sin(x),
		\qquad h \in L^\infty_\lambda\times L^\infty_\lambda,
	$$
	where
	\begin{align*}
		L^{1, \nu}_t(h) & =  -\kappa \int_0^t{ \int_{\mathbb{T} \times \R }
				                                    \left[  \frac{\d}{\d y} \E
					                                    \cos (X^{\delta_y, \omega, f}_{t, \theta})
					                                    \right](D_f v^f_x(\theta, y, \omega) \cdot h)
				                                    \mathcal{L}(X^{\nu^\omega, \omega, f}_\theta) (\d y) g(\d \omega) \d \theta}, \\
		L^{2, \nu}_t(h) & =  -\kappa \int_0^t{ \int_{\mathbb{T} \times \R }  \left[  \frac{\d}{\d y} \E
					                                    \sin (X^{\delta_y, \omega, f}_{t, \theta})
					                                    \right](D_f v^f_x(\theta, y, \omega) \cdot h)
				                                    \mathcal{L}(X^{\nu^\omega, \omega, f}_\theta) (\d y) g(\d \omega) \d \theta}.
	\end{align*}
\end{proposition}
\begin{proof}
	Define for $\theta \in [0, t]$ and $y \in \T$:
	\[ \Psi^\omega_{t, \theta}(y) := \left( \frac{\d}{\d y} \E \cos(X^{\delta_y, \omega, f}_{t, \theta}) \right) (  v^{f+h}_x - v^f_x)(\theta, y, \omega). \]
	Using  Proposition~\ref{prop:sensitivity formula SDE} with $g = \cos$, we have
	\[ \E \cos(X^{\nu^\omega,\omega, f+h}_{t, 0}) - \E \cos(X^{\nu^\omega,\omega, f}_{t, 0}) = - \int_0^t{ \E \Psi^\omega_{t, \theta}(X^{\nu^\omega, \omega, f+h}_{\theta, 0}) \d \theta}. \]
	Applying Proposition~\ref{prop:sensitivity formula SDE} with $g = \Psi^\omega_{t, \theta}$, we find that:
	\begin{align*}
		 & \E \Psi^\omega_{t, \theta}(X^{\nu^\omega, \omega, f+h}_{\theta, 0}) - \E \Psi^\omega_{t, \theta}(X^{\nu^\omega, \omega, f}_{\theta, 0})                                                                                 \\
		 & \quad = -\int_0^\theta \int_{\T} \left( \frac{\d }{\d y} \E_y \Psi^\omega_{t, \theta}(X^{\delta_y, \omega, f}_{\theta, u}) \right) (v_x^{f+h} - v_x^{f})(u, y, \omega) \cL(X^{\nu^\omega, \omega, f+h}_{u})(\d y) \d u.
	\end{align*}
	Let
	\[ B^\omega_{t, \theta} := \sup_{y \in \T} \int_0^\theta e^{\lambda(\theta-u)}   \left| \frac{\d }{\d y} \E_y \Psi^\omega_{t, \theta}(X^{\delta_y, \omega, f}_{\theta, u}) \right| \d u.
	\]
	We deduce that:
	\[ \left|  \E \Psi^\omega_{t, \theta}(X^{\nu, \omega, f+h}_{\theta, 0}) - \E \Psi^\omega_{t, \theta}(X^{\nu, \omega, f}_{\theta, 0}) \right| \leq C B^\omega_{t, \theta} e^{-\lambda \theta} \norm{h}_\lambda.    \]
	By Lemma~\ref{lem:grad estimate}, we obtain that:
	\[  \left| \frac{\d }{\d y} \E_y \Psi^\omega_{t, \theta}(X^{\delta_y, \omega, f}_{\theta, u}) \right|   \leq C \norm{h}_\lambda e^{-\lambda \theta} \left( e^{-\frac{\sigma^2}{2} (t-\theta)} +e^{-\lambda t} \right),   \]
	and so
	\[ B^\omega_{t,\theta} \leq C(\lambda) \norm{h}_\lambda \left( e^{-\frac{\sigma^2}{2}(t-\theta)} + e^{-\lambda t} \right).    \]
	Altogether, we deduce that:
	\[  \left| \int_0^t{ \E \Psi^\omega_{t, \theta}(X^{\nu^\omega, \omega, f+h}_{\theta, 0}) \d \theta} - \int_0^t{ \E \Psi^\omega_{t, \theta}(X^{\nu^\omega, \omega, f}_{\theta, 0}) \d \theta} \right| \leq C(\lambda) \norm{h}_\lambda^2 e^{-\lambda t}. \]
	Consider now:
	\[ \tilde{\Psi}^\omega_{t, \theta}(y) := \left( \frac{\d}{\d y} \E \cos(X^{\delta_y, \omega, f}_{t, \theta}) \right) Dv^{f}_x (\theta, y, \omega) \cdot h, \]
	where $ Dv^{f}_x (\theta, y, \omega) \cdot h$ is given by Lemma~\ref{lem:diff-vf}.
	Using \eqref{eq:estimate vx f+h - vx f - Dvxf}, we deduce that:
	\begin{align*}
		\left| \E \cos(X^{\nu^\omega, \omega, f+h}_t) - \E \cos(X^{\nu^\omega, \omega, f}_t) +
		\int_0^t{\E \tilde{\Psi}_{t, \theta}(X^{\nu^\omega, \omega, f}_\theta) \d \theta} \right|  \leq C e^{-\lambda t} \norm{h}^2_\lambda.
	\end{align*}
	Integrating this estimate with respect to $g(\d \omega)$, this concludes the proof for $\cos$.  The same
	argument is used for $\sin$.
\end{proof}

We next specialize to the $g$-uniform distribution.
\begin{corollary}
	\label{cor:derivative equil Phi}
	For $ h \in L^\infty_\lambda\times L^\infty_\lambda$,
	$$
		(D_f \Phi(\unig, 0) \cdot h) (t, x)
		= -h(t,x) +   (L^{1, \unig} h)(t) \cos(x) + (L^{2, \unig} h)(t) \sin(x),
	$$
	where the linear operators $L^{1, \unig}, L^{2, \unig}:
		L^\infty_\lambda \times L^\infty_\lambda \rightarrow L^\infty_\lambda$ are given by:
	\begin{align*}
		(L^{1, \unig} h) (t) & = \frac{\kappa}{2} \int_0^t \int_\theta^\infty
		e^{-\frac{\sigma^2}{2}(t-\theta)} e^{-\gamma(u-\theta)} \int g(\d \omega)
		\left[ h^1(u) \cos(\omega(t-u)) - h^2(u) \sin (\omega(t-u))  \right] \d u \d \theta, \\
		(L^{2, \unig} h) (t) & = \frac{\kappa}{2} \int_0^t \int_\theta^\infty
		e^{-\frac{\sigma^2}{2}(t-\theta)} e^{-\gamma(u-\theta)} \int g(\d \omega)
		\left[ -h^1(u) \sin(\omega(t-u)) + h^2(u) \cos (\omega(t-u))  \right] \d u \d \theta.
	\end{align*}
\end{corollary}
\begin{remark}
	When $g$ is symmetric around zero, the above formula simplifies and:
	\[     (D_f \Phi(\unig, 0) \cdot h) (t, x)
		= -h(t,x) +   \kappa (L h^1)(t) \cos(x) + \kappa (L h^2)(t) \sin(x),   \]
	where $L: L^\infty_\lambda \rightarrow L^\infty_\lambda$ is given by \eqref{eq:operator-L}.
\end{remark}
\begin{proof}
	We show that the formula above holds for $L^{1, \unig}$, the second formula is shown similarly. First, combining Proposition~\ref{prop:Phi is Fr\'echet differentiable} and Corollary~\ref{cor:derivative-v_x_near_0}, we find that:
	\begin{align*}  L^{1, \unig}_t(h) = & - \kappa \int_0^t \int_{\T \times \R}  \sin(y + \omega(t-\theta) ) e^{-\frac{\sigma^2}{2}(t-\theta)}                                                                     \\
	                                    & \left[ (- K^2_\omega * h^1 + K^1_\omega * h^2)(\theta) \cos(y) - (K^1_\omega * h^1 + K^2_\omega * h^2)(\theta) \sin(y) \right] \frac{\d y}{2 \pi} g(\d \omega) \d \theta\end{align*}
	As $\sin(y + \omega(t-\theta)) = \cos(y) \sin(\omega(t-\theta)) + \sin(y) \cos(\omega(t-\theta))$, we find that:
	\begin{align*}
		L^{1, \unig}_t(h) = & - \frac{\kappa}{2} \int_0^t \int_{\R}  \sin(\omega(t-\theta) ) e^{-\frac{\sigma^2}{2}(t-\theta)}  (- K^2_\omega * h^1 + K^1_\omega * h^2)(\theta)  g(\d \omega) \d \theta     \\
		                    & \quad + \frac{\kappa}{2} \int_0^t \int_{\R}  \cos(\omega(t-\theta) ) e^{-\frac{\sigma^2}{2}(t-\theta)}  (K^1_\omega * h^1 + K^2_\omega * h^2)(\theta)  g(\d \omega) \d \theta \\
		                    & =: A(h^1) + B(h^2).
	\end{align*}
	We focus on $A(h^1)$ and recall that
	$$
		-( K^2_\omega * h^1)(\theta)  = - \int_\theta^\infty e^{-\gamma (u-\theta) } \sin(\omega(u-\theta)) h^1(u) \d u,
		\qquad
		( K^1_\omega * h^1)(\theta)  = \int_\theta^\infty e^{-\gamma (u-\theta) } \cos(\omega(u-\theta)) h^1(u) \d u.
	$$
	We now use the elementary identities,
	\[
		\sin(\omega(t-\theta)) \sin(\omega(u-\theta)) = \frac{1}{2} \cos(\omega(t-u)) - \frac{1}{2} \cos(\omega(t-2\theta +u)),
	\]
	\[ \cos(\omega(t-\theta)) \cos(\omega(u-\theta)) =  \frac{1}{2} \cos(\omega(t-u)) + \frac{1}{2} \cos(\omega(t-2\theta +u))
	\]
	to arrive at
	\[ A(h^1) =  \frac{\kappa}{2} \int_0^t \int_\theta^\infty e^{-\frac{\sigma^2}{2}(t-\theta)} e^{-\gamma(u-\theta)} h^1(u) \int_{\R} g(\d \omega) \cos(\omega(t-u)) \d u \d \theta. \]
	The term $B(h^2)$ is computed similarly, giving the stated formula.
\end{proof}

%%%%%%%%%%%%%%%%%%%%%%%%%%%
\appendix

%%%%%%%%%%%%%%%%%%%%%%%%%%%%%%
\section{Coupled PDEs and the potential structure}
\label{app:pde}

In this appendix, we discuss alternative definitions
of the Kuramoto game given by nonlinear
partial differential equations and show that
it is a potential game by constructing the corresponding
control problem.  Although these could be highly relevant
to some studies, we do not utilize them in our analysis
and state them in an appendix for possible future references.

\subsection{Coupled HJB-KFP system}
\label{subsec:HJB}

The Kuramoto MFG is equivalently
defined through the solutions
of coupled \emph{Hamilton-Jacobi-Bellman} (HJB)
and \emph{Kolmogorov-Fokker-Planck} (KFP) equations
for the value function $v$
and the density $\mu$, with
$v, \mu :[0,\infty) \times \T \times \R \mapsto \R$
given by,
\begin{align*}
	 & \beta v(t,x,\omega) - \partial_t v(t, x,\omega) -
	\frac{\sigma^2}{2} \partial^2_{xx}v(t, x,\omega) -
	\omega \partial_x v(t, x,\omega) +
	\frac{1}{2} (\partial_x v(t, x,\omega))^2
	= \kappa c(x, \mu(t,\cdot,\cdot)),                   \\
	 & \partial_t \mu(t, x,\omega) -
	\frac{\sigma^2}{2} \partial^2_{xx} \mu(t, x,\omega) +
	\partial_x ((\omega -\partial_x v(t, x,\omega)) \mu(t, x,\omega)) = 0,
\end{align*}
where
$$
	c(x, \mu(t,\cdot,\cdot))= c(x, \mu_t),
	\qquad
	\mu_t(\d x, \d \omega)
	= \mu(t,x,\omega)\ \d x\  g(\d \omega),
$$
and $c(x, \mu_t)$ is defined in \eqref{fo:developed_c}.
When the above system has classical solutions,
the flow of probability
measures $\mu_t \in \cP(\T \times \R)$
defined above is the corresponding Nash equilibrium.
Since the above equations do not involve
any derivatives with respect to the $\omega$ variable,
one may view them as parametrized by the frequency $\omega$
with solutions $v^\omega,\mu^\omega : [0,\infty)\times \T \mapsto \R$
defined by $(v^\omega,\mu^\omega)(t,x):= (v,\mu)(t,x,\omega)$.
The coupling in the frequency variable
is through the cost $c(x, \mu_t) = 2 \kappa \int_{\T}  \sin^2((x-y)/2) \mu^\omega_t(\d y) g(\d \omega)$ which uses all frequencies. This is exactly
the same mechanism discussed in Remark~\ref{rem:coupling0}.

\subsection{Nash-Lasry-Lions equation}
\label{subsec:NLL}
Another characterization of MFGs is provided
by the so-called master equation.
Following the terminology introduced in~\cite{HSY},
we would like to refer to it as the \emph{Nash-Lasry-Lions} (NLL)
equation, and to state it, we set
$\cX := \T \times \R$.
Then, according to \cite{CD,HSY},
the tagged player's optimal value function $v : [0,\infty) \times  \cX \times \cP(\cX)
	\mapsto \R$ formally satisfies,
\begin{align*}
	\beta v(t,x,\omega,\mu)
	- \partial_ t v(t,x,\omega,\mu) = &
	\frac{\sigma^2}{2} \partial^2_{xx} v(t, x,\omega,\mu) +
	\omega \partial_x v(t, x,\omega,\mu) -
	\frac{1}{2} (\partial_x v(t, x,\omega,\mu))^2                                                            \\
	                                  & + \int_\cX \cI(v)(t,x,\omega,\mu)(x',\omega')\ \mu(d x', \d \omega')
	+  \kappa c(x,\mu),
\end{align*}
where
$$
	\cI(v)(t,x,\omega,\mu)(x',\omega')=
	\frac{\sigma^2}{2}  D_{x'}\partial_\mu v(t,x,\omega,\mu)(x',\omega')
	+ (\omega' - v_x(t,x',\omega',\mu))
	\partial_\mu v(t,x,\omega,\mu)(x',\omega'),
$$
the operator $\partial_\mu $ is the Lions-derivative in the measure
variable $\mu$ which
results in a  function on $(x',\mu')\in \cX$
(see \cite{CD} for its definition and properties), and $D_{x'}$ is the
derivative with respect to the $x'$-variable.

Since all measures on $\cP(\cX)$ that we consider have
second marginals equal to $g \in \cP(\R)$, this subset
of $\cP(\cX)$ is invariant for this equation.
Using this property, we may consider functions
on the set $[0,\infty)\times \cX \times \cM$,
where $ \cM:= \cM(\R \mapsto \cP(\T))$ is the set of
measurable functions of frequencies to probability
measures on $\T$.  For $m \in \cM$, $\omega \in \R$,
we have $m(\omega)\in \cP(\T)$.
Then, in this variable, the integral term in the equation
takes the form
$$
	\int_\cX \cI(v)(t,x,\omega,m)(x',\omega') \ m(\omega')(d x')\, g(\d \omega').
$$

\subsection{Potential structure}
\label{subsec:potential}
It is known that the problem with homogeneous frequencies
is a potential game. In the simple setting
of \cite{FS}, it is explicitly stated and studied.
This structure holds in our model as well.
Indeed, consider the
McKean-Vlasov optimal control
problem of minimizing
$$
	J(\alpha) := \int_0^\infty
	e^{-\beta t} (\ \frac12 \E[\alpha_t^2]
	+ \kappa \Phi(\mu_t^\alpha)\ ) \ \d t
$$
over all adapted processes $\alpha$,
where $\mu^\alpha_t \in \cP(\T \times \R)$ is
the law of the couple $(X_t^{\alpha,\omega},\omega)$ and
$$
	\d X^{\alpha,\omega}_t = (\omega +\alpha_t) \d t + \sigma \d B_t.
$$
Additionally, the intrinsic frequency $\omega$
is $\cF_0$-measurable and has distribution $g$,
the initial condition $\mu_0 =\cL(X_0,\omega)$ satisfies
$\mu_0(\T,\d \omega)=g(\d \omega)$, and
for any $\nu \in \cP(\T\times\R)$
$$
	\Phi(\nu)= \int_\T  \int_\T \sin^2(\frac12[x-x'])\
	\nu(\d x, \R) \ \nu(\d x', \R).
$$
The second marginal of
$\mu_t^\alpha$
does not change in time and is always equal to $g$:
$\mu_t^\alpha(\T,\d \omega)=g(\d \omega)$ for all $t$.

\begin{lemma}[Theorem 3.6 \cite{HS1}]
	\label{lem:potential}
	The minimizers of the above McKean-Vlasov optimal control
	problem are Nash equilibria of the
	Kuramoto MFG with intrinsic frequencies.
\end{lemma}
This follows from the fact that
$$
	\delta_\nu\Phi(\nu)(x,\omega)= c(x,\nu),
	\qquad \forall\ (x,\omega)\in \cX,
$$
where $\delta_\nu\Phi$ is the linear derivative of $\Phi$.

\begin{remark}
	\label{rem:coupling}
	The coupling between the intrinsic frequencies
	is only through the fact that the interaction present in the cost function relies on
	the joint law $\mu_t^\alpha$ of the couple $(X_t^{\alpha,\omega},\omega)$.
	The same mechanism is also at play in the coupling between
	the HJB and KFP equations introduced in the previous subsection,
	as well as in the original definition as discussed in Remark~\ref{rem:coupling0}.
\end{remark}

%\vspace{-2em}
%%%%%%%%%%%%%%%%%%%%%%%%%%%
\section{The implicit function theorem}
We recall the version of the implicit theorem we use.
{
\begin{theorem}[Implicit function theorem]
	\label{th:IFT}
	Let $\Lambda$ be a set equipped with a gauge function $\delta: \Lambda \rightarrow [0, \infty)$ and let $\lambda_0 \in \Lambda$ be such that $\delta(\lambda_0) = 0$. Let $(X, \|\cdot\|)$ be a Banach space. Let $A \subset \Lambda$ be a subset containing $\lambda_0$ and $B \subset X$ be an open neighborhood of $0$. Let $\Phi: A \times B \rightarrow X$ be such that $\Phi(\lambda_0, 0) = 0$ and:
	\begin{enumerate}\itemsep=-1pt
		\item $\lim_{\delta(\lambda) \to 0} \|\Phi(\lambda, 0)\| = 0$.
		\item For all $\lambda \in A$, $\Phi(\lambda, \cdot)$ is Fr\'echet differentiable.
		\item $\lim_{\delta(\lambda) \to 0, \, \|x\| \to 0} \|D_x \Phi(\lambda, x) - D_x \Phi(\lambda_0, 0)\| = 0$.
		\item The linear map $D_x \Phi(\lambda_0, 0)$ is invertible with a bounded inverse.
	\end{enumerate}
	Then there exist $\rho, \rho' > 0$ defining an open neighborhood $A_1 = \{\lambda \in A \mid \delta(\lambda) < \rho\}$ and $B_1 = \{x \in B \mid \|x\| < \rho'\}$, and a function $\phi: A_1 \mapsto B_1$ such that
	$$
		\forall (\lambda, x) \in A_1 \times B_1, \quad \Phi(\lambda, x) = 0 \iff x = \phi(\lambda).
	$$
	In addition, if there is a constant $C > 0$ such that $\|\Phi(\lambda, x) - \Phi(\lambda_0, x)\| \leq C \delta(\lambda)$ for all $\lambda \in A_1$ and $x \in B_1$, then there exists a constant $\tilde{C} > 0$ such that $\|\phi(\lambda)\| \leq \tilde{C} \delta(\lambda)$ for all $\lambda \in A_1$.
\end{theorem}
}
\begin{proof}
	We choose $A_1 = \{\lambda \in A \mid \delta(\lambda) < \rho\}$ and $B_1 = \{x \in B \mid \|x\| < \rho'\}$ for small $\rho, \rho'$ to be chosen later. Fix $\lambda \in A_1$ and define $T_\lambda: B_1 \rightarrow X$:
	$$
		T_\lambda(x) = x - \left[ D_x \Phi(\lambda_0, 0) \right]^{-1} \cdot \Phi(\lambda, x).
	$$
	By Assumption 2, $T_\lambda$ is Fr\'echet differentiable with
	$$
		D_x T_\lambda(x) = I_X - \left[ D_x \Phi(\lambda_0, 0) \right]^{-1} \cdot D_x \Phi(\lambda, x).
	$$
	By Assumption 3, we can find $\rho, \rho' > 0$ small enough such that
	$$
		\forall (\lambda,x) \in A_1 \times B_1, \quad \|D_x T_\lambda(x)\| \leq \frac{1}{2}.
	$$
	In addition, using Assumption 1 and Assumption 4, $\|T_\lambda(0)\| \leq \rho' / 2$ for $\rho$ small enough.
	Altogether,
	$$
		\|T_\lambda(x)\| \leq \|T_\lambda(x) - T_\lambda(0)\| + \|T_\lambda(0)\| \leq \frac{1}{2} \|x\| + \frac{\rho'}{2} \leq \rho'.
	$$
	This shows that $T_\lambda$ maps $\bar{B}_1$ into itself. In addition, $T_\lambda$ is $1/2$-Lipschitz on $\bar{B}_1$. We deduce that $T_\lambda$ has a unique fixed point in $\bar{B}_1$, which we denote by $\phi(\lambda)$. Since $\Phi(\lambda_0, 0) = 0$, we easily see that $T_{\lambda_0}(0) = 0$, meaning $\phi(\lambda_0) = 0$.

	For $\lambda \in A_1$, we have:
	\begin{align*}
		\|\phi(\lambda)\| & = \|\phi(\lambda) - \phi(\lambda_0)\|                                                                                                 \\
		                  & = \|T_\lambda(\phi(\lambda)) - T_{\lambda_0}(\phi(\lambda_0))\|                                                                       \\
		                  & \leq \|T_\lambda(\phi(\lambda)) - T_{\lambda_0}(\phi(\lambda))\| + \|T_{\lambda_0}(\phi(\lambda)) - T_{\lambda_0}(\phi(\lambda_0))\|.
	\end{align*}
	As $\left[ D_x \Phi(\lambda_0, 0) \right]^{-1}$ is a bounded linear operator and $\|\Phi(\lambda, x) - \Phi(\lambda_0, x)\| \leq C \delta(\lambda)$, there is a constant $C_1$ such that $\|T_\lambda(\phi(\lambda)) - T_{\lambda_0}(\phi(\lambda))\| \leq C_1 \delta(\lambda)$. In addition, as $T_{\lambda_0}$ is $1/2$-Lipschitz, we deduce that:
	$$
		\|\phi(\lambda)\| \leq C_1 \delta(\lambda) + \frac{1}{2} \|\phi(\lambda) - \phi(\lambda_0)\| = C_1 \delta(\lambda) + \frac{1}{2} \|\phi(\lambda)\|.
	$$
	Subtracting $\frac{1}{2} \|\phi(\lambda)\|$ from both sides yields $\frac{1}{2}\|\phi(\lambda)\| \leq C_1 \delta(\lambda)$, so $\|\phi(\lambda)\| \leq 2 C_1 \delta(\lambda)$. This concludes the proof with $\tilde{C} = 2 C_1$.
\end{proof}
\section{Parabolic estimates}
\label{sec:appendix:parabolic}
In this section, we provide a proof of Proposition~\ref{prop:regular solution HJB}.
To simplify the notation, we assume that $\sigma=1$.
\subsubsection*{Study of the linear non-homogeneous heat equation}
If $v: \mathbb{T} \rightarrow \mathbb{R}$ is an $\alpha$-Hölder continuous function, we write
\[
	\norm{v}_{C^\alpha} := \sup_{x \in \mathbb{T}}{|v(x)|} + \sup_{x,y \in \mathbb{T}, x \neq y} \frac{ |v(x) - v(y)|}{
		|x-y|^\alpha}.
\]
We first recall some parabolic regularity result of the heat equation.
\begin{lemma}
	\label{lem:construction solution heat equation}
	Let $g: \mathbb{R}_+ \times \mathbb{T} \rightarrow \mathbb{R}$ be a bounded measurable function, and for $t\ge
		0$ and $x\in\mathbb{T}$ define
	\[
		(Sg)(t, x) :=  \int_t^\infty{ e^{-\beta (u-t) }  \E g(u, x + \omega (u-t) + B_{u-t}) \d u}.
	\]
	Then for each $\theta \in (0,1)$, there exists a constant $C(\theta)$ such that:
	\begin{enumerate}\itemsep=-1pt
		\item For all $t \geq 0$, $(Sg)(t, \cdot)$ is differentiable, $(Sg)_x(t, \cdot)$ is $\theta$-Hölder and:
		      \begin{equation} \sup_{t \geq 0} \norm{(Sg)_x(t, \cdot)}_{C^\theta} \leq C(\theta) \norm{g}_\infty.
			      \label{eq:theta holder (Sg)x}
		      \end{equation}
		\item  For all $x \in \mathbb{T}$, $(Sg)(\cdot, x)$ is $\theta$-Hölder and for all $x \in \mathbb{T}$ and  $t,t' \in
			      \mathbb{R}_+$,
		      \begin{equation}
			      \label{eq:theta holder (Sg)}
			      |t-t'| \leq 1 \implies |(Sg)(t, x)-(Sg)(t',x)| \leq C(\theta) \norm{g}_\infty  \left( |t-t'|^{\theta}  + \abs{\omega} |t-t'| \right).
		      \end{equation}
		\item For all $x \in \mathbb{T}$, $(Sg)_x(\cdot, x)$ is $\theta/2$-Hölder and for all $x \in \mathbb{T}$ and $t,t' \in
			      \mathbb{R}_+$,
		      \begin{equation}
			      \label{eq:theta holder (Sg)x in t}
			      |t-t'| \leq 1 \implies |(Sg)_x(t, x)-(Sg)_x(t',x)| \leq C(\theta) \norm{g}_\infty(1+\abs{\omega}) |t-t'|^{\theta/2}.
		      \end{equation}
		      The constant $C(\theta)$ is independent of $\omega$.
	\end{enumerate}
\end{lemma}

\begin{proof}
	Let $J(x, u, t) := \E g(u, x + \omega(u-t) + B_{u-t})$. We have
	\[
		J(x, u, t)= \int_{\mathbb{R}}{ g(u, z) G\bigl(\frac{z-x - \omega (u-t)}{\sqrt{u-t}}\bigr)\frac{\d z}{\sqrt{u-t}} }, \quad \text{ with } G(r) :=
		\frac{e^{-r^2/2}}{\sqrt{2 \pi}}.
	\]
	This shows that $J(\cdot, u, t)$ is $C^\infty$ on $u > t$ and:
	\begin{align*}
		J_x(x, u, t) & = -\int_{\mathbb{R}}{ g(u, z)
			                  \frac{1}{\sqrt{u-t}}G'\bigl(\frac{z-x - \omega(u-t)}{\sqrt{u-t}}\bigr)\frac{\d z}{\sqrt{u-t}}} \\
		             & = - \int_{\mathbb{R}}{ \frac{g(u, x + \omega(u-t) + \sqrt{u-t} y)}{\sqrt{u-t}} G'(y) \d y}.
	\end{align*}
	We deduce that $x \mapsto (Sg)(t, x)$ is differentiable with:
	\begin{equation} (Sg)_x(t, x) =  -\int_t^\infty{ \int_{\mathbb{R}}{ e^{-\beta(u-t)} \frac{g(u, x + \omega(u-t) +  \sqrt{u-t}
					y)}{\sqrt{u-t}} G'(y) \d y}\d u}.
		\label{eq:Dx S(g)}
	\end{equation}
	Let $\theta \in (0, 1)$. Because $u \mapsto {(u-t)^{-1/2 - \theta/2}}$ is integrable near $u = t$, we deduce that
	$(Sg)_x$ is $\theta$-Hölder and \eqref{eq:theta holder (Sg)x} holds.
	We now prove that $(Sg)(\cdot, x)$ is $\theta$-Hölder for all $\theta \in (0, 1)$.
	First, for $u > t$, $t \mapsto J(x,u, t)$ is differentiable and there exists a constant $C$ (independent of $\omega$) such that:
	\[ \left| J_t(x, u, t) \right| \leq C \norm{g}_\infty \left( \frac{1}{u-t} + \frac{\abs{\omega}}{\sqrt{u-t}} \right). \]
	Let $\epsilon \in (0, 1)$. We have for some constant $C$:
	\begin{align*}
		|(Sg)(t+\epsilon, x) - (Sg)(t, x)| \leq C \norm{g}_\infty \epsilon + \int_{t+\epsilon}^\infty{e^{-\beta(u-t-\epsilon)} |J(x, u,
			                                                                                       t+\epsilon) - J(x, u, t)| \d u }.
	\end{align*}
	Using Hölder's inequality, we have for $p = 1/(1-\theta)$:
	\begin{align*} | J(x, u, t+\epsilon) - J(x, u, t)| & \leq C \norm{g}_\infty \epsilon^\theta \left[ \int_t^{t+\epsilon}{ (u-s)^{-p}
			                                                                                            \d s} \right]^{1/p} + \frac{C \epsilon \norm{g}_\infty \abs{\omega}}{\sqrt{u-t}}               \\
	                                                   & \leq C \norm{g}_\infty \epsilon^\theta \frac{[2(p-1)]^{1/p}}{(u-(t+\epsilon))^{\theta}} +  \frac{C \epsilon \norm{g}_\infty \abs{\omega}}{\sqrt{u-t}}.
	\end{align*}
	As the right-hand side is integrable with respect to $u$, we deduce that
	\eqref{eq:theta holder (Sg)} holds. The third point is proved similarly: we have
	\[ |J_{xt}(x, u, t)| \leq C \norm{g}_\infty \left( \frac{1}{(u-t)^{3/2}} + \frac{\abs{\omega}}{u-t} \right);  \]
	and the Hölder inequality gives similarly \eqref{eq:theta holder (Sg)x in t}.
\end{proof}
With a mild regularity assumption on $g$, $(Sg)$ is a classical solution of the following heat equation:
\begin{lemma}
	Assume that for some $\alpha \in (0, 1)$, $g$ is $\alpha$-Hölder with respect to $x$, uniformly in $t$. Then
	$(Sg)_x(t, \cdot)$ and $(Sg)(\cdot, x)$ are both differentiable everywhere, and:
	\begin{equation}
		\label{eq:sol edp g}
		(Sg)_t + \frac{1}{2} (Sg)_{xx}  + \omega (Sg)_x+ g = \beta (Sg)
	\end{equation}
	In addition, $(Sg)_{xx}$ is $\alpha$-Hölder with respect to $x$ uniformly in $t$, and $\alpha/2$-Hölder with
	respect to $t$, uniformly in $x$. The associated Hölder constants only depend on $\alpha, \beta, \omega$ and
	$\norm{g}_{C^\alpha}$.
\end{lemma}

\begin{proof}
	We have for all $u > t$:
	\begin{align*}
		J_{xx}(x, u, t)
		 & = \int_{\mathbb{R}} \bigl(g(u, z) - g(u, x + \omega (u-t))\bigr)  \frac{1}{u-t} G''\bigl(\frac{z-x - \omega(u-t)}{\sqrt{u-t}}\bigr)\frac{\d z}{\sqrt{u-t}} \\
		 & \quad + g(u, x + \omega (u-t)) \int_{\mathbb{R}}{ \frac{1}{u-t}G''\bigl(\frac{z-x - \omega(u-t)}{\sqrt{u-t}}\bigr)\frac{\d z}{\sqrt{u-t}}}                 \\
		 & =: A + B.
	\end{align*}
	It holds that $B = 0$. Using that $g$ is $\alpha$-Hölder, we have
	\begin{align*}
		|A|
		 & \leq \int_{\mathbb{R}}{ \frac{|z-x - \omega (u-t)|^\alpha}{ |u-t|^{\alpha/2} |u-t|^{1-\alpha/2}}
			        |G''|\bigl(\frac{z-x - \omega(u-t)}{\sqrt{u-t}}\bigr)\frac{\d z}{\sqrt{u-t}}}            \\
		 & = \int_{\mathbb{R}}{ \frac{y^{\alpha}}{ (u-t)^{1 - \alpha/2}} |G''|(y) \d y} \leq \frac{C_\alpha}{(u-t)^{1 -
				                                                                                     \alpha/2}}.
	\end{align*}
	This last quantity being integrable with respect to $u$ near $u = t$, we deduce that $(Sg)_x(t,
		\cdot)$ is
	differentiable with a
	bounded derivative given by
	\begin{equation}
		\label{eq:representation_of_Sxx}
		(Sg)_{xx}(t, x) = \int_t^\infty{ \int_{\mathbb{R}}{ e^{-\beta(u-t)} \frac{g(u, x + \omega(u-t) + \sqrt{u-t} y)-g(u, x)}{u-t} G''(y)
				\d y}\d u}.
	\end{equation}
	In addition, by \cite[Ch. 4, estimate (2.8)]{MR0241822}, it holds that $(Sg)_{xx}$ is $\alpha$-Hölder uniformly
	in $t$, and so there exists a constant $C$ (only depending on $\alpha$, $\beta$ and $\omega$) such that
	\[ \sup_{t \geq 0} \norm{(Sg)_{xx}(t, \cdot)}_{C^\alpha} \leq C \norm{g}_\alpha.  \]
	Moreover, by  \cite[Ch. 4, estimate (2.9)]{MR0241822}, $(Sg)_{xx}(x, \cdot)$ is $\alpha/2$ Hölder, uniformly
	in $x \in \mathbb{T}$.
	Finally, note that $J(x, u, t)$ solves the heat equation on $\{t < u\} \times \mathbb{T}$:
	\[ J_t(x, u,t) + \omega J_x(x, u, t)+  \frac{1}{2} J_{xx} (x, u, t) = 0. \]
	So it holds that $|J_t(x,u,t)| \leq \frac{\abs{\omega}}{\sqrt{u-t}} + \frac{C_\alpha}{|u-t|^{1-\alpha/2}}$ and we deduce that
	$(Sg)(\cdot, x)$
	is differentiable everywhere, with a derivative given by \eqref{eq:sol edp g}. This concludes the proof.
\end{proof}

%%%%%%%%%%%%%%%%%%%%%%%%%%%%%%%%%%%%%%%%%%%%%%%%%	
\subsubsection*{Construction of a solution using the Schauder fixed point theorem}

Let $\mathcal{Y}$ be the space of all continuous bounded functions $C(\mathbb{R}_+ \times \mathbb{T})$ such
that for all
$t \geq 0$, $u(t, \cdot)$ is $C^1(\mathbb{T})$ and $\sup_{t \geq 0} \sup_{x \in \mathbb{T}} |u_x(t, x)| < \infty$. We equip $\mathcal{Y}$ with the norm
\[
	\norm{u}_{\mathcal{Y}} := \sup_{t \geq 0} \bigl( \norm{u(t, \cdot)}_\infty + \norm{u_x(t, \cdot)}_\infty  \bigr).
\]
This makes $(\mathcal{Y}, \norm{\cdot}_{\mathcal{Y}})$ a Banach space.
Let $\lambda > 0$, $c := \sqrt{\frac{\beta}{2}}$, and $\mathfrak{S}$ be the following closed convex
set of $\mathcal{\mathcal{Y}}$:
\[ \mathfrak{S} := \{ u \in \mathcal{Y}~|~ \forall t \geq 0, \quad    \sup_{x \in \mathbb{T}} |u_x(t,x)| \leq c e^{-\lambda
			t}  \}. \]
Given $f \in \mathcal{X}_\lambda$ satisfying $\norm{f}_\lambda \leq  \beta/4$ and $v \in \mathfrak{S}$, we define
\[
	(Tv)(t,x) := \E \int_t^\infty{ e^{-\beta (u-t) } \bigl[ -\frac{1}{2} v^2_x(u, x + \omega (u-t) + B_{u-t}) + f(u, x + \omega(u-t) + B_{u-t})  \bigr]
		\d u}.
\]
\begin{lemma}
	The function $T$ maps $\mathfrak{S}$ into itself.
\end{lemma}
\begin{proof}
	Given $v \in \mathfrak{S}$, let $g := - \frac{1}{2} v^2_x + f$. It holds that
	\[ |g(t, x)| \leq (\frac{1}{2} c^2 + \beta/4) e^{-\lambda t}. \]
	By Lemma~\ref{lem:construction solution heat equation}, $T(v) \in \mathcal{Y}$ with $(Tv)_x$ given by \eqref{eq:Dx
		S(g)}. In particular we have
	\begin{align*}
		|(Tv)_x(t, x)| & \leq (\frac{1}{2} c^2 + \beta/4) e^{-\lambda t} \int_t^\infty{
			                                                                        \frac{e^{-\beta(u-t)}}{\sqrt{u-t}}\int_{\mathbb{R}}{ |G'(y)|\d y \d u } } \\
		               & = \sqrt{\frac{2}{\beta}} \left(\frac{1}{2} c^2 +
		\beta/4 \right) e^{-\lambda t} = c e^{-\lambda t}.
	\end{align*}
\end{proof}
We use the following extension of the Arzelà-Ascoli theorem on $\mathbb{R}_+$:
\begin{lemma}
	\label{lem:ascoli in R+}
	Let $C, \lambda > 0$ be fixed.
	Let $f_n$ be a sequence of functions in $C(\mathbb{R}_+ \times \mathbb{T})$ such that $(f_n)$ is
	equicontinuous and
	\[ \forall n \in \mathbb{N}, \forall t \geq 0, \forall x \in \mathbb{T}, \quad |f_n(t, x)| \leq C e^{-\lambda t}. \]
	Then, there exists $g \in C([0, \infty) \times \mathbb{T})$ and a sub-sequence $n_i$ such that $\lim_{i
			\rightarrow
			\infty} ||f_{n_i}-g||_\infty = 0$.
\end{lemma}
\begin{proof}
	By the Arzelà–Ascoli theorem, we can find a subsequence $\sigma_1(i)$ such that $f_{\sigma_1(i)}$
	converges
	towards $g_1$ uniformly on $[0,1] \times \mathbb{T}$. By induction, we can find sub-sequences
	$\sigma_k(i)$ of $\sigma_{k-1}(i)$ such that $f_{\sigma_k(i)} \rightarrow_i g_k$, uniformly on $[0,k] \times
		\mathbb{T}$. This construction ensures that for
	any integer $m < k$, $x \in \mathbb{T}$
	\[ \forall t \in [0, m], \quad g_k(t, x) = g_m(t, x). \]
	We denote by $g(t, x)$ the common limit, and we define the diagonal sequence by $\sigma(i) := \sigma_i(i)$.
	The sequence $f_{\sigma(i)}$ is a Cauchy sequence for
	the uniform norm on $\mathbb{R}_+ \times \mathbb{T}$. Indeed, let $\epsilon > 0$,  there exists $A > 0$
	large
	enough such that $2 C e^{-\lambda A} < \epsilon$. So for every $t > A$, we have
	\[|f_{\sigma(i)}(t, x) - f_{\sigma(j)}(t, x)| \leq 2 C e^{-\lambda t} \leq \epsilon. \]
	In addition, the sequence $f_{\sigma(i)}$ converges uniformly towards $g$ on $[0, A] \times \mathbb{T}$,
	and so
	there exists $N$ large enough such that for all $i,j \leq N$, for all $t \in [0, A], x \in \mathbb{T}$, we have
	$|f_{\sigma(i)}(t,x) - f_{\sigma(j)}(t,x)| \leq \epsilon$.
	This concludes the proof.
\end{proof}
\begin{proposition}
	$T$ has a fixed point in $\mathfrak{S}$.
\end{proposition}
\begin{proof}
	First, $T$ is continuous. Let us now verify that $T \mathfrak{S}$ is precompact.
	Let $\theta \in (0,1)$ be fixed. Using \eqref{eq:theta holder (Sg)x}, we deduce that there
	exists a constant $C(\beta, \theta)$ such that for all $v \in \mathfrak{S}$:
	\[ \forall t \geq 0, \quad ||(Tv)(t, \cdot)||_{C^\theta} \leq C(\beta, \theta) e^{-\lambda t}.   \]
	In addition, by \eqref{eq:theta holder (Sg)}, there exists another constant $C(\beta, \theta, \omega)$ such that for all $v \in
		\mathfrak{S}$:
	\[ \forall t,t' \in \mathbb{R}_+, \quad |t-t'| \leq 1 \implies |(Tv)(t, x) - (Tv)(t', x)| \leq C(\beta, \theta, \omega) |t-t'|^\theta. \]
	Using \eqref{eq:theta holder (Sg)x in t}, a similar inequality holds (with $\theta/2$ instead of $\theta$) for $(Tv)_x$.
	Using Lemma~\ref{lem:ascoli in R+} twice, we deduce that $T \mathfrak{S}$ is precompact. So the Schauder fixed
	point
	theorem applies \cite[Cor. 11.2]{MR737190} and $T$ has a fixed point in $\mathfrak{S}$.
\end{proof}
We conclude the proof of Proposition~\ref{prop:regular solution HJB} by Lemma~\ref{lem:construction solution
	heat equation} once more. To prove that $v(t, \cdot) \in C^3(\T)$, we finally use \eqref{eq:representation_of_Sxx}.
\bibliographystyle{abbrvnat}
\bibliography{biblio}
\end{document}